\documentclass[12pt]{amsart}
\usepackage{pgfplots}
\usepackage{tikz}

\usetikzlibrary{automata}
\usetikzlibrary{positioning}
\usepackage{amsmath, amssymb} 
\usepackage{amsthm}
\usepackage{graphicx}
\usepackage[shortlabels]{enumitem}

\usepackage{amsfonts,geometry,hyperref}

\newtheorem{remark}{Remark}[section]

\newtheorem{theorem}{Theorem}[section]

\newtheorem{prob}[theorem]{Problem}
\newtheorem{lemma}[theorem]{Lemma}

\DeclareMathOperator{\Spec}{Spec}
\DeclareMathOperator{\L-Spec}{L-spec}
\DeclareMathOperator{\Q-Spec}{Q-spec}
\DeclareMathOperator{\CN-Spec}{CN-spec}
\DeclareMathOperator{\con}{con}

\allowdisplaybreaks
\begin{document}
\title[Various spectra and energies of subgroup generating bipartite graph]{Various spectra and energies of subgroup generating bipartite graph}

\author[S. Das, A. Erfanian and R. K. Nath]{Shrabani Das, Ahmad Erfanian and Rajat Kanti Nath$^*$}

\address{S. Das, Department of Mathematical Sciences, Tezpur University, Napaam-784028, Sonitpur, Assam, India.
\newline
Department of Mathematics, Sibsagar University, Joysagar-785665, Sibsagar, Assam, India.}

\email{shrabanidas904@gmail.com}

\address{A. Erfanian, Department of Pure Mathematics, Ferdowsi University of Mashhad,
	P.O. Box 1159-91775,
	Mashhad, Iran}
\email{erfanian@um.ac.ir}

\address{R. K. Nath, Department of Mathematical Sciences, Tezpur University, Napaam-784028, Sonitpur, Assam, India.} 
\email{ rajatkantinath@yahoo.com}
\thanks{$^*$Corresponding Author}
\begin{abstract}
Let $L(G)$ be the set of all subgroups of a group $G$. The subgroup generating bipartite graph $\mathcal{B}(G)$ defined on $G$ is a bipartite graph  whose vertex set is partitioned into two sets $G \times G$ and $L(G)$, and two vertices $(a, b) \in G \times G$ and $H \in L(G)$ are adjacent if   $H$ is generated by $a$ and $b$. In this paper, we compute various spectra and  energies of $\mathcal{B}(G)$ and determine whether $\mathcal{B}(G)$ is hypoenergetic, hyperenergetic, CN-hyperenergetic, L-hyperenergetic or Q-hyperenergetic if $G$ is a dihedral group of order $2p$ and $2p^2$ and dicyclic group of order $4p$ and $4p^2$, where $p$ is any prime.
We also show that $\mathcal{B}(G)$ satisfies E-LE conjecture for these groups.
\end{abstract}

\thanks{ }

\maketitle

\section{Introduction}
Let $\mathcal{G}$ be a simple graph with vertex set $V(\mathcal{G})=\{v_i:i=1,2, \ldots, n\}$ and $m$ edges. Let $A(\mathcal{G})$ and $D(\mathcal{G})$ denote the adjacency matrix and degree matrix of $\mathcal{G}$ respectively. The spectrum of $\mathcal{G}$ is the set of eigenvalues of $A(\mathcal{G})$ along with their multiplicities. The Laplacian matrix and signless Laplacian matrix of $\mathcal{G}$ are given by $L(\mathcal{G}):=D(\mathcal{G})-A(\mathcal{G})$ and $Q(\mathcal{G}):=D(\mathcal{G})+A(\mathcal{G})$ respectively. The set of eigenvalues of $L(\mathcal{G})$ and $D(\mathcal{G})$ along with their multiplicities are defined as the Laplacian spectrum (L-spectrum) and signless Laplacian spectrum (Q-spectrum) of $\mathcal{G}$ respectively. The common neighborhood of two distinct vertices $v_i$ and $v_j$, denoted by $C(v_i, v_j)$, is the set of all vertices other than $v_i$ and $v_j$ which are adjacent to both $v_i$ and $v_j$. The common neighborhood matrix of $\mathcal{G}$, denoted by $CN(\mathcal{G})$, is defined as
\[
(CN(\mathcal{G}))_{i, j}= \begin{cases}
	|C(v_i, v_j)|, & \text{if } i \neq j \\
	0, & \text{if } i=j.
\end{cases}
\]
The set of all eigenvalues of $CN(\mathcal{G})$ along with their multiplicities is called the common neighborhood spectrum (CN-spectrum) of $\mathcal{G}$. We write $\Spec(\mathcal{G})$, $\L-Spec(\mathcal{G})$, $\Q-Spec(\mathcal{G})$ and $\CN-Spec(\mathcal{G})$ to denote the spectrum, L-spectrum, Q-spectrum and CN-spectrum of $\mathcal{G}$ respectively. These multi-sets $\Spec(\mathcal{G})$/$\L-Spec(\mathcal{G})$/$\Q-Spec(\mathcal{G})$/$\CN-Spec(\mathcal{G})$ are described as  $\{(x_1)^{n_1}, (x_2)^{n_2}, \dots, (x_k)^{n_k}\}$ where $x_i$'s are eigenvalues of $A(\mathcal{G})$/$L(\mathcal{G})$/$Q(\mathcal{G})$/$CN(\mathcal{G})$ and $n_i$'s are their multiplicities. The graph $\mathcal{G}$ is called integral/L-integral/Q-integral/CN-integral if  $\Spec(\mathcal{G})$/$\L-Spec(\mathcal{G})$/$\Q-Spec(\mathcal{G})$/$\CN-Spec(\mathcal{G})$ contain only integers. In graph theory, it is a general problem to determine integral, L-integral, Q-integral and CN-integral graphs (see \cite{HaSc-1974,BaCv-2002,FKMN-2005,Merris-1994,CRS-2007,Stanic-2009,ASG}).

The energy, $E(\mathcal{G})$ and common neighborhood energy (CN-energy), $E_{CN}(\mathcal{G})$ of $\mathcal{G}$ are the sum of the absolute values of the eigenvalues of $A(\mathcal{G})$ and $CN(\mathcal{G})$ respectively. Thus 
\[
E(\mathcal{G}) = \sum_{\alpha \in \Spec(\mathcal{G})}|\alpha| \quad \text{ and } \quad E_{CN}(\mathcal{G}) = \sum_{\beta \in \CN-Spec(\mathcal{G})}|\beta|.
\]
 The Laplacian energy (L-energy), $LE(\mathcal{G})$ and signless Laplacian energy (Q-energy), $LE^{+}(\mathcal{G})$ of $\mathcal{G}$ are defined as
\[
LE(\mathcal{G})=\sum_{\lambda \in \L-Spec(\mathcal{G})} \left| \lambda-\frac{2m}{n} \right| \quad \text{and} \quad LE^{+}(\mathcal{G})= \sum_{\mu \in \Q-Spec(\mathcal{G})} \left| \mu-\frac{2m}{n} \right|.
\]
Note that energy of a graph was introduced by Gutman \cite{Gutman-78} and   the other graph energies mentioned above were introduced by Gutman et al. \cite{GZ06,ACGMR11,ASG}.
It is well-known that $E(K_n)=LE(K_n)=LE^+(K_n)=2(n-1)$ and $E_{CN}(K_n)=2(n-1)(n-2)$.
A graph $\mathcal{G}$ with $|V(\mathcal{G})| =n$  is called hyperenergetic (see \cite{Gutman-99,Walikar-99}) if $E(\mathcal{G}) > E(K_n)$. It is called hypoenergetic if $E(\mathcal{G}) < n$. Similarly, $\mathcal{G}$ is called L-hyperenergetic (see \cite{fns20}) if $LE(\mathcal{G}) > LE(K_n)$, Q-hyperenergetic (see \cite{fns20}) if $LE^+(\mathcal{G}) > LE^+(K_n)$ and CN-hyperenergetic (see \cite{ASG}) if $E_{CN}(\mathcal{G}) > E_{CN}(K_n)$. 
It is still an open problem to find a CN-hyperenergetic graph (see \cite[Open Problem 1]{ASG}).

The common neighborhood graph \cite{AAGS12} of $\mathcal{G}$, denoted by  $\con(\mathcal{G})$,  is a graph whose adjacency matrix is given by
\begin{align*}
	(A(\con(\mathcal{G})))_{i,j}=\begin{cases}
		1, & \text{if } |C(v_i, v_j)| \geq 1 \text{ and } i \neq j \\
		0, & \text{otherwise.}
	\end{cases}
\end{align*}
It is easy to see that $CN(\mathcal{G}) = A(\con(\mathcal{G}))$  if   $\mathcal{G}=K_{1, n}$, the star on $n+1$ vertices. In Example 2.1 \cite{AAGS12}, it was shown that $\con(K_{1, n}) =  K_1 \sqcup K_n$. Therefore,  $\CN-Spec(K_{1, n})= \Spec(K_1) \cup \Spec(K_n)$.

Since the past few decades, graphs have been defined on finite groups and studied towards the characterization of finite groups/graphs. Some well-researched examples of such graphs include commuting graph, \quad non-commuting graph, \quad power graph, B superA graph etc. (see \cite{cameron2021graphs,CaJa-2024,KSCC-2021,ANC-2022}).
Note that none of these graphs are bipartite in nature. In \cite{DEN-23}, we have introduced a bipartite graph called subgroup generating bipartite graph (in short SGB-graph) as given below. Let $G$ be a finite group and $L(G) := \{H : H \text{ is a subgroup of } G\}$. 
The SGB-graph of $G$, denoted by $\mathcal{B}(G)$, is a graph  whose vertex set $V(\mathcal{B}(G))$ is partitioned into two sets $G \times G$ and $L(G)$, and two vertices $(a, b) \in G \times G$ and $H \in L(G)$ are adjacent if  $H$ is generated by $a$ and $b$. 
In \cite{DEN-23}, we have obtained relations between $\mathcal{B}(G)$  and various probabilities (such as probability generating  a given subgroup, commuting probability, cyclicity degree, nilpotency degree, solvability degree) associated to finite groups. We have established connections between $\mathcal{B}(G)$ and generating graph of $G$. Various graph parameters of $\mathcal{B}(G)$ such as independence number, domination number, girth, diameter etc. were also discussed. In \cite{DEN-24}, various topological indices of $\mathcal{B}(G)$ such as first and second Zagreb indices, Randic connectivity index, Atom-Bond-Connectivity index, Geometric Arithmetic index  etc. were computed and Hansen-Vuki{\v{c}}evi{\'c}  conjecture \cite{hansen2007comparing} was verified. Another type of bipartite graph defined on $G$ has been considered recently by Mahtabi et al. \cite{MEM-2023}.

In this paper, we compute spectrum, L-spectrum,  Q-spectrum, CN-spectrum and their corresponding energies for $\mathcal{B}(G)$ if  $G = D_{2p}, D_{2p^2}, Q_{4p}$ and $Q_{4p^2}$ for any prime $p$, where
$D_{2m}$ is the dihedral group presented by $\langle a, b: a^m=b^2=1, bab=a^{-1} \rangle$ and $Q_{4m}$ is the dicyclic group presented by $\langle a, b : a^{2m} = 1, b^2 = a^m, bab^{-1} = a^{-1} \rangle$. As a consequence, we show that $\mathcal{B}(G)$ is not integral but L-integral, Q-integral and CN-integral if   $G$ is one of the above mentioned groups. Further, it is shown that 
$\mathcal{B}(G)$ is hypoenergetic but neither hyperenergetic, L-hyperenergetic, Q-hyperenergetic nor CN-hyperenergetic if $G$ is one of the above mentioned groups. 
Gutman et al. \cite{E-LE-Gutman}  conjectured that $E(\mathcal{G}) \leq LE(\mathcal{G})$  which is known as E-LE conjecture. Gutman \cite{Gutman-78} also conjectured that ``$\mathcal{G}$ is not hyperenergetic if $\mathcal{G} \ncong K_{|v(\mathcal{G})|}$". However, both the conjectures were disproved. We show that $\mathcal{B}(G)$ satisfies both the above mentioned conjectures if $G = D_{2p}, D_{2p^2}, Q_{4p}$ and $Q_{4p^2}$. We conclude the paper showing that  $\mathcal{B}(G)$ satisfies E-LE conjecture if $G$ is one of the above mentioned graphs.  The motivation of this work comes from the study of spectral aspects of commuting/non-commuting/commuting conjugacy class/power and super commuting graphs of groups. For instance, various spectrum and energies of commuting graphs were computed in \cite{AL2021,DN-2017,DN-2017MV,N-2018,TA19,DBN-2020,DN-2021,DN-2022,FSN-2021,NFDS-2021}; non-commuting graphs in \cite{DN18,DDN18,SN24,FN24,AEN17,GGA17,GGAZ17}, commuting conjugacy class graphs in \cite{JN25,BN24,BN21},  power graphs in \cite{MGA17,HA17,P19,RPCA23,CPA18,S22,BA23} and super commuting graphs in \cite{DMP24,ACGN24}.

In Section 2, we recall the graph structures of $\mathcal{B}(G)$ when $G = D_{2p}, D_{2p^2}, Q_{4p}$ and $Q_{4p^2}$, and various spectra and energies of a star graph. In section 3, we consider spectral aspects of $\mathcal{B}(D_{2p})$ and $\mathcal{B}(D_{2p^2})$ while in Section 4, we consider spectral aspects of $\mathcal{B}(Q_{4p})$ and $\mathcal{B}(Q_{4p^2})$.

\section{Some useful results}
We write $mK_{1, r}$ to denote $m$ copies of the star $K_{1, r}$. The following results obtained in \cite{DEN-23}  are useful in our computations.

\begin{theorem}\label{structure_of_D_2p}
\cite[Theorem 4.1-- 4.2]{DEN-23} Let $D_{2n}=\langle a, b: a^n=b^2=1, bab=a^{-1} \rangle$ be the dihedral group of order $2n$ ($n \geq 3$) and $p$ be a prime. Then 
	\[
	\mathcal{B}(D_{2n})= \begin{cases}
		K_2 \sqcup pK_{1, 3} \sqcup K_{1, p^2-1} \sqcup K_{1, 3p(p-1)}, & \text{ when } n = p\\
		K_2 \sqcup p^2K_{1, 3} \sqcup K_{1, p^2-1} \sqcup K_{1, p^4-p^2} \sqcup pK_{1, 3p(p-1)} \sqcup K_{1, 3p^2(p^2-p)}, & \text{ when } n = p^2.
	\end{cases}
	\] 
 \end{theorem}
 \begin{theorem}\label{structure_of_Q_4p}
 \cite[Theorem 4.3]{DEN-23}	Let $Q_{4p} =  \langle a, b : a^{2p} = 1, b^2 = a^p, bab^{-1} = a^{-1} \rangle$ be the dicyclic group of order $4p$, where $p$ is a prime. Then
 	\[
 	\mathcal{B}(Q_{4p})=\begin{cases}
 		K_2 \sqcup K_{1, 3} \sqcup 3K_{1, 12} \sqcup K_{1, 24}, & \text{ when } p=2 \\
 		K_2 \sqcup K_{1, 3} \sqcup pK_{1, 12} \sqcup K_{1, p^2-1} \sqcup K_{1, 3p^2-3} \sqcup K_{1, 12p^2-12p}, & \text{ when } p \geq 3.
 	\end{cases}
 	\]
 \end{theorem}
 \begin{theorem}\label{structure_of_Q_4p^2}
 \cite[Theorem 4.4]{DEN-23}	Let $Q_{4p^2} =  \langle a, b : a^{2p^2} = 1, b^2 = a^{p^2}, bab^{-1} = a^{-1} \rangle$ be the dicyclic group of order $4p^2$, where $p$ is a prime. Then
 	\[
 	\mathcal{B}(Q_{4p^2})=\begin{cases}
 		K_2 \sqcup K_{1, 3} \sqcup 5K_{1, 12} \sqcup 2K_{1, 24} \sqcup K_{1, 48} \sqcup K_{1, 96}, & \text{ when } p=2 \\
 		K_2 \sqcup K_{1, 3} \sqcup p^2K_{1, 12} \sqcup K_{1, p^2-1} \sqcup K_{1, 3p^2-3} \sqcup K_{1, 3p^4-3p^2} \\ \qquad \qquad \sqcup (p-1)K_{1, 12p^2-12p} \sqcup K_{1, 13p^4-12p^3+11p^2-12p}, & \text{ when } p \geq 3.
 	\end{cases}
 	\]
 \end{theorem}
The following well-known lemmas which give various spectra and energies of a star are also useful in computing various spectra and energies of $\mathcal{B}(G)$.
\begin{lemma}\label{all_spec_of_star}
	If \,$\mathcal{G}=K_{1, n}$, the star on $n+1$ vertices, then
	$\Spec(\mathcal{G})=\big{\{}(0)^{n-1}, \left(\sqrt{n}\right)^1,$ $ \left(-\sqrt{n}\right)^1\big{\}}$, 
	$\L-Spec(\mathcal{G})=\left\lbrace(0)^1, (1)^{n-1}, (n+1)^1\right\rbrace = \Q-Spec(\mathcal{G})$ and
	$\CN-Spec(\mathcal{G})= \big{\{}(0)^1, $ $(-1)^{n-1}, (n-1)^1\big{\}}$.
\end{lemma}
\begin{lemma}\label{all_energy_of_star}
	If \,$\mathcal{G}=K_{1, n}$, the star on $n+1$ vertices, then
	$E(\mathcal{G})=2 \sqrt{n}$, 
	$LE(\mathcal{G})=\frac{2n^2+2}{n+1} = LE^+(\mathcal{G})$
	and 
	$E_{CN}(\mathcal{G})= 2n-2$. 
\end{lemma}


\section{Certain dihedral groups}
In this section, we compute spectrum, Laplacian spectrum,  signless Laplacian spectrum, common neighborhood spectrum and their corresponding energies of SGB-graph of dihedral groups of order $2p$ and $2p^2$ for any prime $p$. We also determine if $\mathcal{B}(D_{2p})$ and $\mathcal{B}(D_{2p^2})$ are hypoenergetic, hyperenergetic, L-hyperenergetic, Q-hyperenergetic and CN-hyperenergetic. Finally, we show that $\mathcal{B}(D_{2p})$ and $\mathcal{B}(D_{2p^2})$ satisfy E-LE conjecture.
\begin{theorem}\label{all_spec_of_B(D_2p)}
	If $G=D_{2p}$ 
	then 
	\begin{align*}
		\Spec(\mathcal{B}(G))&=\left\lbrace(0)^{4p^2-p-3}, (-1)^1, (1)^1, \left(\sqrt{3}\right)^p, \left(-\sqrt{3}\right)^p, \left(\sqrt{p^2-1}\right)^1, \left(-\sqrt{p^2-1}\right)^1,\right. \\
		&\qquad\left. \left(\sqrt{3p^2-3p}\right)^1, \left(-\sqrt{3p^2-3p}\right)^1\right\rbrace, 
	\end{align*}        
	\[
	\L-Spec(\mathcal{B}(G))=\left\lbrace(0)^{p+3}, (1)^{4p^2-p-3}, (2)^1, (4)^p, (p^2)^1, (3p^2-3p+1)^1\right\rbrace = \Q-Spec(\mathcal{B}(G)) 
	\]
and	$\CN-Spec(\mathcal{B}(G))=\left\lbrace(0)^{p+4}, (-1)^{4p^2-p-3}, (2)^p, (p^2-2)^1, (3p^2-3p-1)^1\right\rbrace.$
\end{theorem}
\begin{proof}
	By Theorem \ref{structure_of_D_2p}, we have $\mathcal{B}(G)=K_2 \sqcup pK_{1, 3} \sqcup K_{1, p^2-1} \sqcup K_{1, 3p(p-1)}$. Again, by Lemma \ref{all_spec_of_star}, we get  
	\[
	\Spec(K_2) = \left\lbrace (-1)^1, (1)^1\right\rbrace, \quad \Spec(pK_{1,3}) =\left\lbrace(0)^{2p}, \left(\sqrt{3}\right)^p, \left(-\sqrt{3}\right)^p\right\rbrace,
	\]
	\[
	\Spec(K_{1,p^2-1})= \left\lbrace(0)^{p^2-2}, \left(\sqrt{p^2-1}\right)^1, \left(-\sqrt{p^2-1}\right)^1\right\rbrace  \text{ and }
	\] 
	\[   
	\Spec(K_{1,3p(p-1)}) = \left\lbrace(0)^{3p^2-3p-1}, \left(\sqrt{3p^2-3p}\right)^1, \left(-\sqrt{3p^2-3p}\right)^1\right\rbrace.
	\]
	Since $\Spec(\mathcal{B}(G)) = \Spec(K_2) \cup \Spec(pK_{1,3}) \cup \Spec(K_{1,p^2-1}) \cup \Spec(K_{1,3p(p-1)})$ we get the required expression for $\Spec(\mathcal{B}(G))$.
	By Lemma \ref{all_spec_of_star}, we also get 
	\[
	\L-Spec(K_2) = \Q-Spec(K_2) =\left\lbrace(0)^1, (2)^1\right\rbrace,
	\]
	\[
	\L-Spec(pK_{1,3}) = \Q-Spec(pK_{1,3}) =\left\lbrace (0)^p, (1)^{2p}, (4)^p\right\rbrace, 
	\]
	\[
	\L-Spec(K_{1,p^2-1})= \Q-Spec(K_{1,p^2-1})= \left\lbrace (0)^1, (1)^{p^2-2}, (p^2)^1\right\rbrace
	\text{ and}
	\]
	\[
	\L-Spec(K_{1,3p(p-1)})= \Q-Spec(K_{1,3p(p-1)})= \left\lbrace(0)^1, (1)^{3p^2-3p-1}, (3p^2-3p+1)^1\right\rbrace.
	\] 
	Since 
	\begin{align*}
		\L-Spec(\mathcal{B}(G))&=\L-Spec(K_2) \cup\L-Spec(pK_{1,3}) \cup \L-Spec(K_{1,p^2-1}) \cup \L-Spec(K_{1,3p(p-1)}) \\
		&= \Q-Spec(K_2) \cup \Q-Spec(pK_{1,3}) \cup \Q-Spec(K_{1,p^2-1}) \cup \Q-Spec(K_{1,3p(p-1)})\\
		&= \Q-Spec(\mathcal{B}(G))
	\end{align*}
	we get the required expressions for $\L-Spec(\mathcal{B}(G))$ and $\Q-Spec(\mathcal{B}(G))$.
		Further,
		\[
		\CN-Spec(K_2) = \left\lbrace (0)^2\right\rbrace, \quad
		\CN-Spec(pK_{1,3}) = \left\lbrace (0)^p, (-1)^{2p}, (2)^p\right\rbrace, 
		\]
		\[
		\CN-Spec(K_{1,p^2-1})=   \left\lbrace (0)^1, (-1)^{p^2-2}, (p^2-2)^1\right\rbrace \quad \text{ and}
		\]
		\[
		\CN-Spec(K_{1,3p(p-1)}) = \left\lbrace (0)^1, (-1)^{3p^2-3p-1}, (3p^2-3p-1)^1\right\rbrace.
		\] 
		Since  $\CN-Spec(\mathcal{B}(G))=\CN-Spec(K_2) \cup \CN-Spec(pK_{1,3}) \cup \CN-Spec(K_{1,p^2-1})$ we get the required expression for $\CN-Spec(\mathcal{B}(G))$.
		\end{proof}
\begin{theorem}\label{all_spec_of_B(D_2p^2)}
	If $G=D_{2p^2}$  
	then 
	\begin{align*}
		\Spec(\mathcal{B}(G))=&\left\lbrace (0)^{4p^4-p^2-p-4}, (-1)^1, (1)^1, \left(\sqrt{3}\right)^{p^2}, \left(-\sqrt{3}\right)^{p^2}, \left(\sqrt{p^2-1}\right)^1, \right. \\
		&\quad\left(-\sqrt{p^2-1}\right)^1, \left(\sqrt{p^4-p^2}\right)^1, \left(-\sqrt{p^4-p^2}\right)^1, \left(\sqrt{3p^2-3p}\right)^p,  \\
		&\qquad \qquad\qquad\left.\left(-\sqrt{3p^2-3p}\right)^p, \left(\sqrt{3p^4-3p^3}\right)^1, \left(-\sqrt{3p^4-3p^3}\right)^1\right\rbrace,
	\end{align*}
		\begin{align*}
		\L-Spec(\mathcal{B}(G))=&\Big{\{}(0)^{p^2+p+4}, (1)^{4p^4-p^2-p-4}, (2)^1, (4)^{p^2}, (p^2)^1, (p^4-p^2+1)^1, \\
		&\qquad \qquad\qquad\quad (3p^2-3p+1)^{p}, (3p^4-3p^3+1)^{1}\Big{\}} = \Q-Spec(\mathcal{B}(G))
	\end{align*}
 and
 \begin{align*} 
		\CN-Spec(\mathcal{B}(G))=&\Big{\{}(0)^{p^2+p+5}, (-1)^{4p^4-p^2-p-4}, (2)^{p^2}, (p^2-2)^1, (p^4-p^2-1)^1, \\
		&\qquad\qquad\qquad\qquad\qquad\qquad  (3p^2-3p-1)^p, (3p^4-3p^3-1)^1\Big{\}}. 
	\end{align*}
\end{theorem}
\begin{proof}
	By Theorem \ref{structure_of_D_2p}, we have $\mathcal{B}(G)=K_2 \sqcup p^2K_{1, 3} \sqcup K_{1, p^2-1} \sqcup K_{1, p^4-p^2} \sqcup pK_{1, 3p(p-1)} \sqcup K_{1, 3p^2(p^2-p)}$.  Again, by Lemma \ref{all_spec_of_star}, we get 
	\[
	\Spec(K_2) =\left\{(-1)^1, (1)^1\right\}, \,  
	\Spec(p^2K_{1, 3}) =  \left\{(0)^{2p^2}, \left(\sqrt{3}\right)^{p^2}, \left(-\sqrt{3}\right)^{p^2}\right\},
	\]
	\[
	\Spec(K_{1, p^2-1}) = \left\{(0)^{p^2-2}, \left(\sqrt{p^2-1}\right)^1, \left(-\sqrt{p^2-1}\right)^1\right\}, \]
	\[ 
	\Spec(K_{1, p^4-p^2}) = \left\{(0)^{p^4-p^2-1}, \left(\sqrt{p^4-p^2}\right)^1, \left(-\sqrt{p^4-p^2}\right)^1\right\},
	\]
	\[
	\Spec(pK_{1, 3p(p-1)}) = \left\{ (0)^{p(3p^2-3p-1)}, \left(\sqrt{3p^2-3p}\right)^p, \left(-\sqrt{3p^2-3p}\right)^p\right\} \quad \text{and}
	\]
	\[
	\Spec(K_{1, 3p^2(p^2-p)}) = \left\{(0)^{3p^4-3p^3-1}, \left(\sqrt{3p^4-3p^3}\right)^1, \left(-\sqrt{3p^4-3p^3}\right)^1\right\}.
	\]
	Since 
	
\noindent	$\Spec(\mathcal{B}(G))=\Spec(K_2) \cup \Spec(p^2K_{1, 3}) \cup \Spec(K_{1, p^2-1}) \cup \Spec(K_{1, p^4-p^2}) \cup \Spec(pK_{1, 3p(p-1)})$ we get the required expression for $\Spec(\mathcal{B}(G))$.
		By Lemma \ref{all_spec_of_star}, we also get 
	\[
	\L-Spec(K_2) = \Q-Spec(K_2)= \left\{(0)^1, (2)^1\right\}, 
	\]
	\[
	 \L-Spec(p^2K_{1, 3}) = \Q-Spec(p^2K_{1, 3}) =\left\{(0)^{p^2},  (1)^{2p^2}, (4)^{p^2}\right\},
	\]
	\[
	\L-Spec(K_{1, p^2-1}) = \Q-Spec(K_{1, p^2-1}) = \left\{(0)^1, (1)^{p^2-2}, (p^2)^1\right\},
	\]
	\[
	\L-Spec(K_{1, p^4-p^2}) = \Q-Spec(K_{1, p^4-p^2}) = \left\{(0)^1, (1)^{p^4-p^2-1}, \left(p^4-p^2+1 \right)^1\right\},
	\]
	\[
	\L-Spec(pK_{1, 3p(p-1)}) = \Q-Spec(pK_{1, 3p(p-1)}) = \left\{(0)^p, (1)^{p(3p^2-3p-1)}, \left(3p^2-3p+1\right)^p\right\} \, \text{and} 
	\]
	\[
	\L-Spec(K_{1, 3p^2(p^2-p)}) = \Q-Spec(K_{1, 3p^2(p^2-p)}) = \left\{ (0)^1, (1)^{3p^4-3p^3-1}, \left(3p^4-3p^3+1\right)^1\right\}.
	\]
	Since
	\begin{align*}
		\L-Spec(\mathcal{B}(G))= &\L-Spec(K_2) \cup \L-Spec(p^2K_{1, 3}) \cup \L-Spec(K_{1, p^2-1}) \cup \L-Spec(K_{1, p^4-p^2}) \\ 
		& ~~~~~~ \cup \L-Spec(pK_{1, 3p(p-1)}) \cup \L-Spec(K_{1, 3p^2(p^2-p)}) \\
		= & \Q-Spec(K_2) \cup \Q-Spec(p^2K_{1, 3}) \cup \Q-Spec(K_{1, p^2-1}) \cup \Q-Spec(K_{1, p^4-p^2}) \\ 
		& ~~~~~~ \cup \Q-Spec(pK_{1, 3p(p-1)}) \cup \Q-Spec(K_{1, 3p^2(p^2-p)})  \\
		= & \Q-Spec(\mathcal{B}(G))  
	\end{align*}
	we get the required expressions for  	$\L-Spec(\mathcal{B}(G))$ and $\Q-Spec(\mathcal{B}(G))$.
	Further,
	\[
	\CN-Spec(K_2) =  \left\{(0)^2\right\}, \quad \CN-Spec(p^2K_{1, 3}) = \left\{ (0)^{p^2}, (-1)^{2p^2}, (2)^{p^2} \right\}, 
	\]
	\[
	\CN-Spec(K_{1, p^2-1}) = \left\{ (0)^1, (-1)^{p^2-2}, (p^2-2)^1\right\},
	\]
	\[
	\CN-Spec(K_{1, p^4-p^2}) = \left\{(0)^1, (-1)^{p^4-p^2-1}, \left(p^4-p^2-1 \right)^1\right\},
	\]
	\[
	\CN-Spec(pK_{1, 3p(p-1)}) = \left\{(0)^p, (-1)^{p(3p^2-3p-1)}, \left(3p^2-3p-1\right)^p\right\} \, \text{and}
	\]
	\[
	\CN-Spec(K_{1, 3p^2(p^2-p)}) = \left\{(0)^1, (-1)^{3p^4-3p^3-1}, \left(3p^4-3p^3-1\right)^1\right\}.
	\]
	Since \quad $\CN-Spec(\mathcal{B}(G)) \,\, = \,\, \CN-Spec(K_2) \,\, \cup \,\, \CN-Spec(p^2K_{1, 3}) \,\, \cup \,\, \CN-Spec(K_{1, p^2-1}) \,\, \cup$\\ $
	\CN-Spec(K_{1, p^4-p^2})$ $\cup \CN-Spec(pK_{1, 3p(p-1)}) \cup \CN-Spec(K_{1, 3p^2(p^2-p)})$ \quad we get the required expression for $\CN-Spec(\mathcal{B}(G))$.
\end{proof}
\begin{theorem}\label{all_energy_of_B(D_2p)}
	If $G=D_{2p}$ 
	then 
	$E(\mathcal{B}(G))= 2+2p\sqrt{3}+2\sqrt{p^2-1}+2\sqrt{3p(p-1)}$,
	$LE(\mathcal{B}(G))= LE^+(\mathcal{B}(G))=\frac{32p^4+2p^2+12p+18}{4p^2+p+3}$ 
	\quad  and
	$E_{CN}(\mathcal{B}(G))=8p^2-2p-6$.
\end{theorem}
\begin{proof}
	By Lemma \ref{all_energy_of_star} we get
	
\noindent	$E(K_2) = 2,   E( pK_{1, 3}) = 2p\sqrt{3},  E(K_{1, p^2-1}) = 2\sqrt{p^2 - 1} \text{ and }  E(K_{1, 3p(p-1)}) = 2\sqrt{3p(p-1)}.$
	
\noindent	Since $E(\mathcal{B}(G)) = E(K_2)+ E(pK_{1, 3})+ E(K_{1, p^2-1})+ E(K_{1, 3p(p-1)})$ we get the required expression for $E(\mathcal{B}(G))$. By Theorem \ref{all_spec_of_B(D_2p)}, we have $\L-Spec(\mathcal{B}(G))=\left\lbrace(0)^{p+3}, (1)^{4p^2-p-3}, (2)^1, \right.$   $(4)^p, (p^2)^1, (3p^2-3p+1)^1\Big{\}} = \Q-Spec(\mathcal{B}(G))$, $m=4p^2$ and $n=4p^2+p+3$. Now, $\left |0-\frac{2(4p^2)}{4p^2+p+3}\right|=\frac{8p^2}{4p^2+p+3}, \left |1-\frac{2(4p^2)}{4p^2+p+3}\right |=\frac{4p^2-p-3}{4p^2+p+3}, \left |2-\frac{2(4p^2)}{4p^2+p+3}\right |=\frac{2p+6}{4p^2+p+3}, \left |4-\frac{2(4p^2)}{4p^2+p+3}\right |=\frac{8p^2+4p+12}{4p^2+p+3}, \left |p^2-\frac{2(4p^2)}{4p^2+p+3}\right |=\frac{4p^4+p^3-5p^2}{4p^2+p+3}$ and $\left |(3p^2-3p+1)-\frac{2(4p^2)}{4p^2+p+3}\right |=\frac{12p^4-9p^3+2p^2-8p+3}{4p^2+p+3}$. Therefore, by definition of (signless) Laplacian energy we have 
 \begin{align*}
     LE(\mathcal{B}(G))&=LE^+(\mathcal{B}(G))\\
     &=(p+3)\frac{8p^2}{4p^2+p+3}+(4p^2-p-3)\frac{4p^2-p-3}{4p^2+p+3}+\frac{2p+6}{4p^2+p+3} \\
     & \qquad \qquad +p\frac{8p^2+4p+12}{4p^2+p+3}+\frac{4p^4+p^3-5p^2}{4p^2+p+3}+\frac{12p^4-9p^3+2p^2-8p+3}{4p^2+p+3} \\
     &=\frac{32p^4+2p^2+12p+18}{4p^2+p+3}.
 \end{align*}
Finally,
$
E_{CN}(K_2) = 0, \,  E_{CN}( pK_{1, 3}) = 4p, \, E_{CN}(K_{1, p^2-1}) = 2p^2 - 4 \text{ and }  E_{CN}(K_{1, 3p(p-1)}) = 6p^2 - 6p -2.$
	Since $E_{CN}(\mathcal{B}(G)) = E_{CN}(K_2)+ E_{CN}(pK_{1, 3})+ E_{CN}(K_{1, p^2-1})+ E_{CN}(K_{1, 3p(p-1)})$ we get the required expression for $E_{CN}(\mathcal{B}(G))$.
\end{proof}
\begin{theorem}\label{all_energy_of_B(D_2p^2)}
	If $G=D_{2p^2}$  
	then $E(\mathcal{B}(G))= 2+2p^2\sqrt{3}+(2p + 2)\sqrt{p^2-1}+4p\sqrt{3p(p-1)},$ $ LE(\mathcal{B}(G)) = \frac{32p^8+2p^4+4p^3+18p^2+16p+32}{4p^4+p^2+p+4}= LE^+(\mathcal{B}(G))$
	and  $E_{CN}(\mathcal{B}(G))=8p^4-2p^2-2p-8$.
\end{theorem}

\begin{proof}
	By Lemma \ref{all_energy_of_star} we get
\begin{center}
	$
	E(K_2) = 2, \, E(p^2K_{1, 3}) = 2p^2\sqrt{3}, \,  E(K_{1, p^2-1}) = 2\sqrt{p^2 - 1}, \,  E(K_{1, p^4-p^2})  = 2p\sqrt{p^2-1},
$
\end{center}
\begin{center}
	$
	E(pK_{1, 3p(p-1)}) = 2p\sqrt{3p(p-1)} \text{ and }
	E(K_{1, 3p^2(p^2-p)}) = 2p \sqrt{3p(p-1)}.
$
\end{center}
	Since $E(\mathcal{B}(G))=E(K_2) + E(p^2K_{1, 3}) + E(K_{1, p^2-1}) + E(K_{1, p^4-p^2}) + E(pK_{1, 3p(p-1)})  + E(K_{1, 3p^2(p^2-p)})$ we get the required expression for $E(\mathcal{B}(G))$. By Theorem \ref{all_spec_of_B(D_2p^2)}, we have 
	
\noindent	$\L-Spec(\mathcal{B}(G))\!=\!\left\lbrace(0)^{p^2+p+4}\!, (1)^{4p^4-p^2-p-4}, (2)^1, (4)^{p^2}, (p^2)^1, (p^4-p^2+1)^1, (3p^2-3p+1)^{p},\right.$  $(3p^4-3p^3+1)^{1}\Big{\}} = \Q-Spec(\mathcal{B}(G))$, $m=4p^4$ and $n=4p^4+p^2+p+4$. Now, 
	
\noindent	$\left |0-\frac{2(4p^4)}{4p^4+p^2+p+4}\right| =\frac{8p^4}{4p^4+p^2+p+4}$, \quad $\left |1-\frac{2(4p^4)}{4p^4+p^2+p+4}\right |=\frac{4p^4-p^2-p-4}{4p^4+p^2+p+4}$, \quad $\left |2-\frac{2(4p^4)}{4p^4+p^2+p+4}\right |$ $=\frac{2p^2+2p+8}{4p^4+p^2+p+4}$, \quad $\left |4-\frac{2(4p^4)}{4p^4+p^2+p+4}\right|$ $=\frac{8p^4+4p^2+4p+16}{4p^4+p^2+p+4}$, \quad $\left |p^2-\frac{2(4p^4)}{4p^4+p^2+p+4}\right |=\frac{4p^6-7p^4+p^3-4p^2}{4p^4+p^2+p+4}$, 

\noindent $\left |(p^4-p^2+1)-\frac{2(4p^4)}{4p^4+p^2+p+4}\right |=$  $\frac{4p^8-3p^6+p^5-p^4-p^3-3p^2+p+4}{4p^4+p^2+p+4}, \left |(3p^2-3p+1)-\frac{2(4p^4)}{4p^4+p^2+p+4}\right |=\frac{12p^6-12p^5-p^4+10p^2-11p+4}{4p^4+p^2+p+4}$ and $\left |(3p^4-3p^3+1)-\frac{2(4p^4)}{4p^4+p^2+p+4}\right |=\frac{12p^8-12p^7+3p^6+5p^4-12p^3+p^2+p+4}{4p^4+p^2+p+4}$. Therefore, by definition of (signless) Laplacian energy we have 
 \begin{align*}
     LE(\mathcal{B}(G))&=LE^+(\mathcal{B}(G))\\
     &=(p^2+p+4)\frac{8p^4}{4p^4+p^2+p+4}+(4p^4-p^2-p-4)\frac{4p^4-p^2-p-4}{4p^4+p^2+p+4} \\
     & \qquad \quad +\frac{2p^2+2p+8}{4p^4+p^2+p+4}+p^2 \times \frac{8p^4+4p^2+4p+16}{4p^4+p^2+p+4}+\frac{4p^6-7p^4+p^3-4p^2}{4p^4+p^2+p+4}\\
     & \qquad \quad +\frac{4p^8-3p^6+p^5-p^4-p^3-3p^2+p+4}{4p^4+p^2+p+4}\\
     & \qquad \qquad \qquad \qquad +p \times \frac{12p^6-12p^5-p^4+10p^2-11p+4}{4p^4+p^2+p+4}\\
     & \qquad \qquad \qquad \qquad +\frac{12p^8-12p^7+3p^6+5p^4-12p^3+p^2+p+4}{4p^4+p^2+p+4} \\
     &=\frac{32p^8+2p^4+4p^3+18p^2+16p+32}{4p^4+p^2+p+4}.
 \end{align*}
 
	Finally,
	
\noindent 	$E_{CN}(K_2) = 0, \, E_{CN}(p^2K_{1, 3}) = 4p^2,   E_{CN}(K_{1, p^2-1}) = 2p^2 - 4,   E_{CN}(K_{1, p^4-p^2})  = 2p^4- 2p^2 - 2,$
	\[
	E_{CN}(pK_{1, 3p(p-1)}) = 6p^3-6p^2 - 2p\text{ and }
	E_{CN}(K_{1, 3p^2(p^2-p)}) = 6p^4 - 6p^3 - 2.
	\]
	Since \, \, \, $E_{CN}(\mathcal{B}(G))=E_{CN}(K_2) + E_{CN}(p^2K_{1, 3}) + E_{CN}(K_{1, p^2-1}) + E_{CN}(K_{1, p^4-p^2})+ E_{CN}(pK_{1, 3p(p-1)})$ $  + E_{CN}(K_{1, 3p^2(p^2-p)})$ we get the required expression for $E_{CN}(\mathcal{B}(G))$.
\end{proof}

\begin{theorem}\label{all_energy_comparison_D_2p}
	If $G=D_{2p}$,  
	then $\mathcal{B}(G)$ is hypoenergetic but not hyperenergetic,  L-hyperenergetic, Q-hyperenergetic and CN-hyperenergetic.
\end{theorem}
\begin{proof}
	By Theorem \ref{structure_of_D_2p} and Theorem \ref{all_energy_of_B(D_2p)}, we have 
	\[
	|V(\mathcal{B}(G))|=4p^2+p+3 \quad \text{ and} \quad E(\mathcal{B}(G))=2+2p\sqrt{3}+2\sqrt{p^2-1}+2\sqrt{3p^2-3p}.
	\]
	Since $4p>2p\sqrt{3}, \, p^2>2\sqrt{p^2-1}$ and $3p^2-3p+1 >2\sqrt{3p^2-3p}$ \,\, we have 
	\[
	4p+p^2+3p^2-3p+1+2 > 2+2p\sqrt{3}+2\sqrt{p^2-1}+2\sqrt{3p^2-3p}.
	\]
	Therefore,
	\begin{equation} \label{D-2p-hypo}
		|V(\mathcal{B}(G))|> E(\mathcal{B}(G)).
	\end{equation}	 
	Thus, $\mathcal{B}(G)$ is hypoenergetic. 
	
	We have $E(K_{|V(\mathcal{B}(G))|}) = E(K_{4p^2+p+3})= 2(4p^2+p+3-1)=8p^2+2p+4 >4p^2+p+3 = |V(\mathcal{B}(G))|> E(\mathcal{B}(G))$ (using \eqref{D-2p-hypo}). Therefore,  $\mathcal{B}(G)$ is not hyperenergetic.

	Also,     $LE(K_{4p^2+p+3})=LE^+(K_{4p^2+p+3})=8p^2+2p+4$. From Theorem \ref{all_energy_of_B(D_2p)}, we have $LE(\mathcal{B}(G))=LE^+(\mathcal{B}(G))=\frac{32p^4+2p^2+12p+18}{4p^2+p+3}$. Now, 
	\begin{align*}
		LE(K_{4p^2+p+3})-LE(\mathcal{B}(G))&= LE^+(K_{4p^2+p+3})-LE^+(\mathcal{B}(G))\\ &=8p^2+2p+4-\frac{32p^4+2p^2+12p+18}{4p^2+p+3} \\
		&= \frac{16p^3+40p^2-2p-6}{4p^2+p+3} > 0.
	\end{align*}
Therefore, 
	\[
	LE^+(K_{4p^2+p+3}) = LE(K_{4p^2+p+3}) > LE(\mathcal{B}(G))=LE^+(\mathcal{B}(G)).
	\] 
	Hence,  $\mathcal{B}(G)$ is neither L-hyperenergetic nor  Q-hyperenergetic.

	We have  $E_{CN}(K_{4p^2+p+3})=2(4p^2+p+3-1)(4p^2+p+3-2)=(8p^2+2p+4)(4p^2+p+1) > 8p^2-2p-6 = E_{CN}(\mathcal{B}(G))$ (using Theorem \ref{all_energy_of_B(D_2p)}). 
	Hence, $\mathcal{B}(G)$ is not CN-hyperenergetic. This completes the proof.
\end{proof}
\begin{theorem}\label{all_energy_comparison_D_2p^2}
	If $G=D_{2p^2}$ then $\mathcal{B}(G)$ is hypoenergetic but not hyperenergetic,  L-hyperenergetic, Q-hyperenergetic and CN-hyperenergetic.
\end{theorem}
\begin{proof}
	By Theorem \ref{structure_of_D_2p} and Theorem \ref{all_energy_of_B(D_2p^2)}, we have $|V(\mathcal{B}(G))|=4p^4+p^2+p+4$ and
	\begin{align*}
		E(\mathcal{B}(G))& = 2+2p^2\sqrt{3}+(2p+2)\sqrt{p^2-1}+4p\sqrt{3p(p-1)}\\
		& = 2+2p^2\sqrt{3}+2\sqrt{p^2-1}+2\sqrt{p^4-p^2}+p^2\sqrt{3p(p-1)}+2\sqrt{3p^4-3p^3}.
	\end{align*}
	Since $4p^2>2p^2\sqrt{3}, \, p^2>2\sqrt{p^2-1}, \, p^4-p^2+1>2\sqrt{p^4-p^2}, \, p(3p^2-3p+1) > 2p\sqrt{3p^2-3p}$ and $3p^4-3p^3+1 >2\sqrt{3p^4-3p^3}$ \, we have 
	\begin{align*}
		4p^2+p^2+p^4-&p^2+1+p(3p^2-3p+1)+3p^4-3p^3+1+2\\
		&> 2+2p^2\sqrt{3}+2\sqrt{p^2-1}+2\sqrt{p^4-p^2}+p^2\sqrt{3p(p-1)}+2\sqrt{3p^4-3p^3}.
	\end{align*}
	Therefore, $|V(\mathcal{B}(G))|> E(\mathcal{B}(G))$. Thus, $\mathcal{B}(G)$ is hypoenergetic. 
	
	We have $E(K_{4p^4+p^2+p+4})= 2(4p^4+p^2+p+4-1)=8p^4+2p^2+2p+6 >4p^4+p^2+p+4 = |V(\mathcal{B}(G))|> E(\mathcal{B}(G))$ (as shown above).  
	Therefore,  $\mathcal{B}(G)$ is not hyperenergetic.
	
	Also,  $LE(K_{4p^4+p^2+p+4})=LE^+(K_{4p^4+p^2+p+4})=8p^4+2p^2+2p+6$. From Theorem \ref{all_energy_of_B(D_2p^2)} we have $LE(\mathcal{B}(G))=LE^+(\mathcal{B}(G))=\frac{32p^8+2p^4+4p^3+18p^2+16p+32}{4p^4+p^2+p+4}$. Now, 
	\begin{align*}
		&LE(K_{4p^4+p^2+p+4})-LE(\mathcal{B}(G)) =LE^+(K_{4p^4+p^2+p+4})- LE^+(\mathcal{B}(G))\\
		&=8p^4+2p^2+2p+6-\frac{32p^8+2p^4+4p^3+18p^2+16p+32}{4p^4+p^2+p+4}\\
		&=\frac{16p^6+16p^5+56p^4-2p^2-2p-8}{4p^4+p^2+p+4} >0.
	\end{align*}
	Therefore, $LE^+(\mathcal{B}(G)) = LE(\mathcal{B}(G))<LE(K_{4p^4+p^2+p+4}) = LE^+(K_{4p^4+p^2+p+4})$. Thus, $\mathcal{B}(G)$ is not L-hyperenergetic and  Q-hyperenergetic. 
 
  We have $E_{CN}(K_{4p^4+p^2+p+4})=2(4p^4+p^2+p+4-1)(4p^4+p^2+p+4-2)=(8p^4+2p^2+2p+6)(4p^4+p^2+p+2) > 8p^4-2p^2-2p-8 = E_{CN}(\mathcal{B}(G))$ (using Theorem \ref{all_energy_of_B(D_2p^2)}). Hence, $\mathcal{B}(G)$ is not CN-hyperenergetic. This completes the proof.
\end{proof}
We conclude this section with the following comparison.
\begin{theorem}
	If $G=D_{2p}$ or $D_{2p^2}$ then $E(\mathcal{B}(G))<LE(\mathcal{B}(G))=LE^+(\mathcal{B}(G))$.
\end{theorem}
\begin{proof}
	\textbf{Case 1.} $G=D_{2p}$.
	
	By Theorem \ref{all_energy_of_B(D_2p)}, we get
	\[
	E(\mathcal{B}(G)) =2+2p\sqrt{3}+2\sqrt{p^2-1}+2\sqrt{3p(p-1)} \quad \text{and}
	\]
	\[
	LE(\mathcal{B}(G))= LE^+(\mathcal{B}(G)) =\frac{32p^4+2p^2+12p+18}{4p^2+p+3}.
	\] 
	
 Also, from Theorem \ref{all_energy_comparison_D_2p}, we have $|V(\mathcal{B}(G))|=4p^2+p+3> E(\mathcal{B}(G))$. Now, 
 \begin{align*}
     LE(\mathcal{B}(G))-|V(\mathcal{B}(G))|&= LE^+(\mathcal{B}(G))-|V(\mathcal{B}(G))| \\
     &=\frac{32p^4+2p^2+12p+18}{4p^2+p+3}-4p^2+p+3 \\
     &= \frac{16p^4-8p^3-23p^2+6p+9}{4p^2+p+3} \\
     &= \frac{p^2(16p^2-8p-23)+6p+9}{4p^2+p+3}:=\frac{f(p)}{g(p)}.
 \end{align*}
 Since, for all prime $p$, $8p(2p-1)>23$  and $4p^2+p+3>0$ we have $\frac{f(p)}{g(p)}>0$. Hence, $LE(\mathcal{B}(G))= LE^+(\mathcal{B}(G))>|V(\mathcal{B}(G))|>E(\mathcal{B}(G))$.

	
\noindent \textbf{Case 2.} $G=D_{2p^2}$. 
	
	By Theorem \ref{all_energy_of_B(D_2p^2)}, we get
	\[
	E(\mathcal{B}(G)) =2+2p^2\sqrt{3}+(2p + 2)\sqrt{p^2-1}+4p\sqrt{3p(p-1)} \quad \text{and}
	\]
	\[
	LE(\mathcal{B}(G)) = LE^+(\mathcal{B}(G)) =\frac{32p^8+2p^4+4p^3+18p^2+16p+32}{4p^4+p^2+p+4}.
	\]
  Also, from Theorem \ref{all_energy_comparison_D_2p^2}, we have $|V(\mathcal{B}(G))|=4p^4+p^2+p+4> E(\mathcal{B}(G))$. Now, 
 \begin{align*}
     LE(\mathcal{B}(G))-|V(\mathcal{B}(G))|&= LE^+(\mathcal{B}(G))-|V(\mathcal{B}(G))| \\
     &=\frac{32p^8+2p^4+4p^3+18p^2+16p+32}{4p^4+p^2+p+4}-4p^4+p^2+p+4 \\
     &= \frac{16p^8-8p^6-8p^5-31p^4+2p^3+9p^2+8p+16}{4p^4+p^2+p+4} \\
     &= \frac{p^4(16p^4-8p^2-8p-31)+2p^3+9p^2+8p+16}{4p^4+p^2+p+4}:=\frac{f(p)}{g(p)}.
 \end{align*}
Since, for all prime $p$, $8p(2p^3-p-1)>31$ and $4p^4+p^2+p+4>0$ we have $\frac{f(p)}{g(p)}>0$. Hence, $LE(\mathcal{B}(G))= LE^+(\mathcal{B}(G))>|V(\mathcal{B}(G))|>E(\mathcal{B}(G))$.
\end{proof}


\section{Certain dicyclic groups}
In this section, we compute spectrum, Laplacian spectrum,  signless Laplacian spectrum, common neighborhood spectrum and their corresponding energies of SGB-graph of dicyclic groups of order $4p$ and $4p^2$ for any prime $p$. We also determine if $\mathcal{B}(Q_{4p})$ and $\mathcal{B}(Q_{4p^2})$ are hypoenergetic, hyperenergetic, L-hyperenergetic, Q-hyperenergetic and CN-hyperenergetic.
\begin{theorem}\label{all_spec_of_B(Q_4p)}
Let $G=Q_{4p}$ and $p$ be any prime.
\begin{enumerate}
\item If $p=2$ then
\[
        \Spec(\mathcal{B}(G))=\left\{(0)^{58}, (-1)^1, (1)^1, (-\sqrt{3})^1, (\sqrt{3})^1, (-\sqrt{12})^3, (\sqrt{12})^3, (-\sqrt{24})^1, (\sqrt{24})^1 \right\},
\]   
\[ \L-Spec(\mathcal{B}(G))=\left\{(0)^6, (1)^{58}, (2)^1, (4)^1, (13)^3, (25)^1 \right\}=\Q-Spec(\mathcal{B}(G)) \,  \text{ and }
\]
        \[
        \CN-Spec(\mathcal{B}(G))=\left\{ (-1)^{58}, (0)^7, (2)^1, (11)^3, (23)^1 \right\}.
        \]
        \item If $p \geq 3$ then
        \begin{align*}
            \Spec&(\mathcal{B}(G))=\left\{(-1)^1, (0)^{16p^2-p-5}, (1)^1, (-\sqrt{3})^1, (\sqrt{3})^1, (-\sqrt{12})^p, (\sqrt{12})^p, (-\sqrt{p^2-1})^1, \right. \\
            &\left.(\sqrt{p^2-1})^1, (-\sqrt{3p^2-3})^1, (\sqrt{3p^2-3})^1, (-\sqrt{12p^2-12p})^1, (\sqrt{12p^2-12p})^1 \right \},
        \end{align*}
        \begin{align*}
            \L-Spec(\mathcal{B}(G))&= \left\{ (0)^{p+5}, (1)^{16p^2-p-5}, (2)^1, (4)^1, (13)^p, (p^2)^1, (3p^2-2)^1, (12p^2-12p+1)^1 \right\} \\
            &= \Q-Spec(\mathcal{B}(G))
        \end{align*}
        \begin{align*}
 \text{ and }       \CN-Spec(\mathcal{B}(G))&=\Big{\{}(-1)^{16p^2-p-5}, (0)^{p+6}, (2)^1, (11)^p, (p^2-2)^1,\\   
        & \qquad \qquad \qquad \qquad \qquad \qquad \qquad (3p^2-4)^1,  (12p^2-12p-1)^1\Big{\}}.
        \end{align*}
    \end{enumerate}
    \end{theorem}
    \begin{proof}
(a) If $p=2$, by Theorem \ref{structure_of_Q_4p}, we have $\mathcal{B}(G)=K_2 \sqcup K_{1, 3} \sqcup 3K_{1, 12} \sqcup K_{1, 24}$. Again, by Lemma \ref{all_spec_of_star}, we get 
        \[
        \Spec(K_2)=\{(-1)^1, (1)^1\}, \quad \Spec(K_{1, 3})=\{(0)^2, (-\sqrt{3})^1, (\sqrt{3})^1\},
        \]
        \[
        \Spec(3K_{1, 12})=\{(0)^{33}, (-\sqrt{12})^3, (\sqrt{12})^3\} \, \text{ and } \, \Spec(K_{1, 24})=\{(0)^{23}, (-\sqrt{24})^1, (\sqrt{24})^1\}.
        \]
        Since $\Spec(\mathcal{B}(G))=\Spec(K_2) \cup \Spec(K_{1, 3}) \cup \Spec(3K_{1, 12}) \cup \Spec(K_{1, 24})$ we get the required expression for $\Spec(\mathcal{B}(G))$. By Lemmma \ref{all_spec_of_star}, we also get
        \[
        \L-Spec(K_2)=\Q-Spec(K_2)=\{(0)^1, (2)^1\},
        \]
        \[
        \L-Spec(K_{1, 3})=\Q-Spec(K_{1, 3})=\{(0)^1, (1)^2, (4)^1\},
        \]
        \[
        \L-Spec(3K_{1, 12})=\Q-Spec(3K_{1, 12})=\{(0)^3, (1)^{33}, (13)^3\} \qquad \text{ and }
        \]
        \[
        \L-Spec(K_{1, 24})=\Q-Spec(K_{1, 24})=\{(0)^1, (1)^{23}, (25)^1\}.
        \]
        Since 
        \begin{align*}
            \L-Spec(\mathcal{B}(G))&=\L-Spec(K_2) \cup \L-Spec(K_{1, 3}) \cup  \L-Spec(3K_{1, 12}) \cup \L-Spec(K_{1, 24}) \\
            &= \Q-Spec(K_2) \cup \Q-Spec(K_{1, 3}) \cup \Q-Spec(3K_{1, 12}) \cup \Q-Spec(K_{1, 24}) \\
            &= \Q-Spec(\mathcal{B}(G))
        \end{align*}
        we get the required expressions for $\L-Spec(\mathcal{B}(G))$ and $\Q-Spec(\mathcal{B}(G))$. Further,
        \[
        \CN-Spec(K_2)=\{(0)^2\}, \qquad \CN-Spec(K_{1, 3})=\{(0)^1, (-1)^2, (2)^1\},
        \]
        \[
        \CN-Spec(3K_{1, 12})=\{(0)^3, (-1)^{33}, (11)^3\} \, \text{ and } \, \CN-Spec(K_{1, 24})=\{(0)^1, (-1)^{23}, (23)^1\}.
        \]
        Since $\CN-Spec(\mathcal{B}(G))=\CN-Spec(K_2) \cup \CN-Spec(K_{1, 3}) \cup \CN-Spec(3K_{1, 12}) $ \\ $\cup \CN-Spec(K_{1, 24})$ we get the required expression for $\CN-Spec(\mathcal{B}(G))$.
        
 (b)  If $p \geq 3$ then by Theorem \ref{structure_of_Q_4p} we have $\mathcal{B}(G)=K_2 \sqcup K_{1, 3} \sqcup pK_{1, 12} \sqcup K_{1, p^2-1} \sqcup K_{1, 3p^2-3} \sqcup K_{1, 12p^2-12p}$. Again, by Lemma \ref{all_spec_of_star}, we get
       \[
        \Spec(K_2)=\{(-1)^1, (1)^1\}, \quad \Spec(K_{1, 3})=\{(0)^2, (-\sqrt{3})^1, (\sqrt{3})^1\},
        \]
        \[
        \Spec(pK_{1, 12})=\left\{(0)^{11p}, (-\sqrt{12})^p, (\sqrt{12})^p\right\},
        \]
        \[
        \Spec(K_{1, p^2-1})=\left\{(0)^{p^2-2}, (-\sqrt{p^2-1})^1, (\sqrt{p^2-1})^1\right\},
        \]
        \[
        \Spec(K_{1, 3p^2-3})=\left\{(0)^{3p^2-4}, (-\sqrt{3p^2-3})^1, (\sqrt{3p^2-3})^1\right\} \qquad \text{ and }
        \]
        \[
        \Spec(K_{1, 12p^2-12p})=\left\{(0)^{12p^2-12p-1}, (-\sqrt{12p^2-12p})^1, (\sqrt{12p^2-12p})^1\right\}.
        \]
        Since $\Spec(\mathcal{B}(G))=\Spec(K_2) \cup \Spec(K_{1, 3}) \cup \Spec(pK_{1, 12}) \cup \Spec(K_{1, p^2-1})$ \\ $\cup \Spec(K_{1, 3p^2-3}) \cup  \Spec(K_{1, 12p^2-12p})$ we get the required expression for $\Spec(\mathcal{B}(G))$. By Lemma \ref{all_spec_of_star}, we also get 
        \[
        \L-Spec(K_2)=\Q-Spec(K_2)=\{(0)^1, (2)^1\},
        \]
        \[
        \L-Spec(K_{1, 3})=\Q-Spec(K_{1, 3})=\{(0)^1, (1)^2, (4)^1\},
        \]
        \[
        \L-Spec(pK_{1, 12})=\Q-Spec(pK_{1, 12})=\{(0)^p, (1)^{11p}, (13)^p\},
        \]
        \[
        \L-Spec(K_{1, p^2-1})=\Q-Spec(K_{1, p^2-1})=\left\{(0)^1, (1)^{p^2-2}, (p^2)^1\right\},
        \]
        \[
        \L-Spec(K_{1, 3p^2-3})=\Q-Spec(K_{1, 3p^2-3})=\left\{(0)^1, (1)^{3p^2-4}, (3p^2-2)^1\right\} \qquad \text{ and }
        \]
        \[
         \L-Spec(K_{1, 12p^2-12p})=\Q-Spec(K_{1, 12p^2-12p})=\left\{(0)^1, (1)^{12p^2-12p-1}, (12p^2-12p+1)^1\right\}.
        \]
        Since $\L-Spec(\mathcal{B}(G))=\Q-Spec(\mathcal{B}(G))$ is the union of the above sets, we get the required expressions for $\L-Spec(\mathcal{B}(G))$ and $\Q-Spec(\mathcal{B}(G))$. Further,
        \[
        \CN-Spec(K_2)=\{(0)^2\}, \qquad \CN-Spec(K_{1, 3})=\{(0)^1, (-1)^2, (2)^1\},
        \]
        \[
        \CN-Spec(pK_{1, 12})=\{(0)^p, (-1)^{11p}, (11)^p\} , \, \CN-Spec(K_{1, p^2-1})=\{(0)^1, (-1)^{p^2-2}, (p^2-2)^1\},
        \]
        \[
        \CN-Spec(K_{1, 3p^2-3})=\left\{(0)^1, (-1)^{3p^2-4}, (3p^2-4)^1\right\} \qquad \text{ and } 
        \]
        \[
        \CN-Spec(K_{1, 12p^2-12p})=\left\{(0)^1, (-1)^{12p^2-12p-1}, (12p^2-12p-1)^1\right\}.
        \]
        Therefore, the union of the above sets gives us the required expression for  $\CN-Spec(\mathcal{B}(G))$.  
    \end{proof}

    \begin{theorem}\label{all_spec_of_B(Q_4p^2)}
 Let $G=Q_{4p^2}$ and $p$ be any prime.
    \begin{enumerate}
        \item If $p=2$ then
        \begin{align*}
        \Spec(\mathcal{B}(G))&=\left\{(0)^{245}, (-1)^1, (1)^1, (-\sqrt{3})^1, (\sqrt{3})^1, (-\sqrt{12})^5, (\sqrt{12})^5, (-\sqrt{24})^2, (\sqrt{24})^2 \right. \\
        & \qquad \qquad \qquad \qquad \qquad \qquad \qquad \left. (-\sqrt{48})^1, (\sqrt{48})^1, (-\sqrt{96})^1, (\sqrt{96})^1 \right\},
        \end{align*}    
         \[ \L-Spec(\mathcal{B}(G))=\left\{(0)^{11}, (1)^{245}, (2)^1, (4)^1, (13)^5, (25)^2, (49)^1, (97)^1 \right\}=\Q-Spec(\mathcal{B}(G)) \,  \text{ and }
        \]
        \[
        \CN-Spec(\mathcal{B}(G))=\left\{ (-1)^{245}, (0)^{12}, (2)^1, (11)^5, (23)^2, (47)^1, (95)^1 \right\}.
        \]
        \item If $p \geq 3$ then
        \begin{align*}
            \Spec(\mathcal{B}(G))&=\left\{(-1)^1, (0)^{16p^4-p^2-p-5}, (1)^1, (-\sqrt{3})^1, (\sqrt{3})^1, (-\sqrt{12})^p, (\sqrt{12})^p,  \right. \\
            &\left. \quad (-\sqrt{p^2-1})^1, (\sqrt{p^2-1})^1, (-\sqrt{3p^2-3})^1, (\sqrt{3p^2-3})^1, (-\sqrt{3p^4-3p^2})^1,  \right. \\ 
            & \left. \qquad \qquad \qquad \quad (\sqrt{3p^4-3p^2})^1, (-\sqrt{12p^2-12p})^{p-1}, (\sqrt{12p^2-12p})^{p-1}, \right. \\
            & \qquad \left. (-\sqrt{13p^4-12p^3+11p^2-12p})^1, (\sqrt{13p^4-12p^3+11p^2-12p})^1 \right \},
        \end{align*}
        \begin{align*}
                  \L-Spec(\mathcal{B}(G))&= \Big{\{} (0)^{p^2+p+5}, (1)^{16p^4-p^2-p-5}, (2)^1, (4)^1, (13)^{p^2}, (p^2)^1, (3p^2-2)^1,  \\
                  & \quad  (3p^4-3p^2+1)^1, (12p^2-12p+1)^{p-1}, (13p^4-12p^3+11p^2-12p+1)^1 \Big{\}} \\
            &= \Q-Spec(\mathcal{B}(G))
        \end{align*}
        \begin{align*}
  \text{and }       \CN-Spec&(\mathcal{B}(G))=\Big{\{}(-1)^{16p^4-p^2-p-5}, (0)^{p^2+p+6}, (2)^1, (11)^{p^2}, (p^2-2)^1, (3p^2-4)^1,  \\
        &  (3p^4-3p^2-1)^1, (12p^2-12p-1)^{p-1},  (13p^4-12p^3+11p^2-12p-1)^1 \Big{\}}.
        \end{align*}
    \end{enumerate}
    \end{theorem}
    \begin{proof}
(a) If $p=2$, then by Theorem \ref{structure_of_Q_4p^2}, we have $\mathcal{B}(G)=K_2 \sqcup K_{1, 3} \sqcup 5K_{1, 12} \sqcup 2K_{1, 24} \sqcup K_{1, 48} \sqcup K_{1, 96}$. Again, by Lemma \ref{all_spec_of_star}, we get 
            \[
	\Spec(K_2) =\left\{(-1)^1, (1)^1\right\}, \,  
	\Spec(K_{1, 3}) =  \left\{(0)^{2},         
       (\sqrt{3})^{1}, (-\sqrt{3})^{1}\right\},
	\]
        \[
        \Spec(5K_{1, 12}) =  \left\{(0)^{55},         
       (\sqrt{12})^{5}, (-\sqrt{12})^{5}\right\}, \, \Spec(2K_{1, 24}) =  \left\{(0)^{46}, (\sqrt{24})^{2}, (-\sqrt{24})^{2}\right\},
        \]
        \[
        \Spec(K_{1, 48}) =  \left\{(0)^{47}, (\sqrt{48})^{1}, (-\sqrt{48})^{1}\right\} \text{ and } \Spec(K_{1, 96}) =  \left\{(0)^{95}, (\sqrt{96})^{1}, (-\sqrt{96})^{1}\right\}.
        \]
        Since $\Spec(\mathcal{B}(G))$ is the union of the above sets, we get the required expression for $\Spec(\mathcal{B}(G))$. By Lemma \ref{all_spec_of_star}, we also get 
	\[
	\L-Spec(K_2) = \Q-Spec(K_2)= \left\{(0)^1, (2)^1\right\},  
	\]
	\[
	 \L-Spec(K_{1, 3}) = \Q-Spec(K_{1, 3}) =\left\{(0)^{1},  (1)^{2}, (4)^{1}\right\},
	\]
        \[
        \L-Spec(5K_{1, 12}) = \Q-Spec(5K_{1, 12}) =\left\{(0)^{5},  (1)^{55}, (13)^{5}\right\},
        \]
        \[
	 \L-Spec(2K_{1, 24}) = \Q-Spec(2K_{1, 24}) =\left\{(0)^{2},  (1)^{46}, (25)^{2}\right\},
	\]  
        \[
	 \L-Spec(K_{1, 48}) = \Q-Spec(K_{1, 48}) =\left\{(0)^{1},  (1)^{47}, (49)^{1}\right\} \text{ and }
	\]
        \[
	 \L-Spec(K_{1, 96}) = \Q-Spec(K_{1, 96}) =\left\{(0)^{1},  (1)^{95}, (97)^{1}\right\}.
	\]
Since $\L-Spec(\mathcal{B}(G))= \Q-Spec(\mathcal{B}(G))$ is the union of the above sets, we get the required expressions for  	$\L-Spec(\mathcal{B}(G))$ and $\Q-Spec(\mathcal{B}(G))$. Further,
	\[
	\CN-Spec(K_2) =  \left\{(0)^2\right\}, \quad \CN-Spec(K_{1, 3}) = \left\{ (0)^{1}, (-1)^{2}, (2)^{1} \right\}, 
	\]
        \[
        \CN-Spec(5K_{1, 12}) = \left\{ (0)^{5}, (-1)^{55}, (11)^{5} \right\}, \, \CN-Spec(2K_{1, 24}) = \left\{ (0)^{2}, (-1)^{46}, (23)^{2} \right\},
        \]
        \[
        \CN-Spec(K_{1, 48}) = \left\{ (0)^{1}, (-1)^{47}, (47)^{1} \right\} \text{ and } \CN-Spec(K_{1, 96}) = \left\{ (0)^{1}, (-1)^{95}, (95)^{1} \right\}.
        \]
Since $\CN-Spec(\mathcal{B}(G))$ is the union of the above sets, we get the required expression for $\CN-Spec(\mathcal{B}(G))$.

(b) By Theorem \ref{structure_of_Q_4p^2}, we have $\mathcal{B}(G)=K_2 \sqcup K_{1, 3} \sqcup p^2K_{1, 12} \sqcup K_{1, p^2-1} \sqcup K_{1, 3p^2-3} \sqcup K_{1, 3p^4-3p^2} \sqcup (p-1)K_{1, 12p^2-12p} \sqcup K_{1, 13p^4-12p^3+11p^2-12p}$.  Again, by Lemma \ref{all_spec_of_star}, we get 
	\[
	\Spec(K_2) =\left\{(-1)^1, (1)^1\right\}, \,  
	\Spec(K_{1, 3}) =  \left\{(0)^{2}, \left(\sqrt{3}\right)^{1}, \left(-\sqrt{3}\right)^{1}\right\},
	\]
 \[
 \Spec(p^2K_{1, 12}) =  \left\{(0)^{11p^2}, \left(\sqrt{12}\right)^{p^2}, \left(-\sqrt{12}\right)^{p^2}\right\},
 \]
	\[
	\Spec(K_{1, p^2-1}) = \left\{(0)^{p^2-2}, \left(\sqrt{p^2-1}\right)^1, \left(-\sqrt{p^2-1}\right)^1\right\}, \]
 \[
 \Spec(K_{1, 3p^2-3}) =  \left\{(0)^{3p^2-4}, \left(\sqrt{3p^2-3}\right)^{1}, \left(-\sqrt{3p^2-3}\right)^{1}\right\},
 \]
	\[ 
	\Spec(K_{1, 3p^4-3p^2}) = \left\{(0)^{3p^4-3p^2-1}, \left(\sqrt{3p^4-3p^2}\right)^1, \left(-\sqrt{3p^4-3p^2}\right)^1\right\},
	\]

\noindent
$\Spec((p-1)K_{1, 12p^2-12p})\!=\!\left\{ (0)^{(p-1)(12p^2-12p-1)}, \left(\sqrt{12p^2-12p}\right)^{p-1}\!, \left(-\sqrt{12p^2-12p}\right)^{p-1}\right\}$

 and
	\begin{align*}
	   \Spec(K_{1, 13p^4-12p^3+11p^2-12p})& = \left\{(0)^{13p^4-12p^3+11p^2-12p-1}, \left(\sqrt{13p^4-12p^3+11p^2-12p}\right)^1, \right. \\
    & \qquad \qquad \qquad \qquad \qquad \left. \left(-\sqrt{13p^4-12p^3+11p^2-12p}\right)^1\right\}. 
	\end{align*}
	
	Since $\Spec(\mathcal{B}(G))$ is the union of the above sets we get the required expression for $\Spec(\mathcal{B}(G))$.
		By Lemma \ref{all_spec_of_star}, we also get 
	\[
	\L-Spec(K_2) = \Q-Spec(K_2)= \left\{(0)^1, (2)^1\right\}, 
	\]
	\[
	 \L-Spec(K_{1, 3}) = \Q-Spec(K_{1, 3}) =\left\{(0)^{1},  (1)^{2}, (4)^{1}\right\},
	\]
        \[
	 \L-Spec(p^2K_{1, 12}) = \Q-Spec(p^2K_{1, 12}) =\left\{(0)^{p^2},  (1)^{11p^2}, (13)^{p^2}\right\},
	\]
	\[
	\L-Spec(K_{1, p^2-1}) = \Q-Spec(K_{1, p^2-1}) = \left\{(0)^1, (1)^{p^2-2}, (p^2)^1\right\},
	\]
        \[
	\L-Spec(K_{1, 3p^2-3}) = \Q-Spec(K_{1, 3p^2-3}) = \left\{(0)^1, (1)^{3p^2-4}, (3p^2-2)^1\right\},
	\]
        \[
	\L-Spec(K_{1, 3p^4-3p^2}) = \Q-Spec(K_{1, 3p^4-3p^2}) = \left\{(0)^1, (1)^{3p^4-3p^2-1}, \left(3p^4-3p^2+1 \right)^1\right\},
	\]
	\begin{align*}
	\L-Spec((p-1)K_{1, 12p^2-12p}) &= \Q-Spec(K_{1, 12p^2-12p}) \\
        &= \left\{(0)^{p-1}, (1)^{(p-1)(12p^2-12p-1)}, \left(12p^2-12p+1\right)^{p-1}\right\} \, \text{ and}     
	\end{align*}
	\begin{align*}
	& \L-Spec(K_{1, 13p^4-12p^3+11p^2-12p}) = \Q-Spec(K_{1, 13p^4-12p^3+11p^2-12p})\\
 & \qquad \qquad \qquad \qquad  = \left\{ (0)^1, (1)^{13p^4-12p^3+11p^2-12p-1}, \left(13p^4-12p^3+11p^2-12p+1\right)^1\right\}.    
	\end{align*}
	Since $\L-Spec(\mathcal{B}(G))= \Q-Spec(\mathcal{B}(G))$ is the union of the above sets, we get the required expressions for  	$\L-Spec(\mathcal{B}(G))$ and $\Q-Spec(\mathcal{B}(G))$.
	Further,
	\[
        \CN-Spec(K_2)=\{(0)^2\}, \qquad \CN-Spec(K_{1, 3})=\{(0)^1, (-1)^2, (2)^1\},
        \]
\noindent
        $\CN-Spec(p^2K_{1, 12})\!=\!\{(0)^{p^2}, (-1)^{11p^2}, (11)^p\},  \CN-Spec(K_{1, p^2-1})\!=\!\{(0)^1, (-1)^{p^2-2}, (p^2-2)^1\},$
        \[
        \CN-Spec(K_{1, 3p^2-3})=\left\{(0)^1, (-1)^{3p^2-4}, (3p^2-4)^1\right\},
        \]
        \[
        \CN-Spec(K_{1, 3p^4-3p^2})=\left\{(0)^1, (-1)^{3p^4-3p^2-1}, (3p^4-3p^2-1)^1\right\},
        \]
        \[
        \CN-Spec((p-1)K_{1, 12p^2-12p})=\left\{(0)^{p-1}, (-1)^{(p-1)(12p^2-12p-1)}, (12p^2-12p-1)^{p-1}\right\} \text{ and }
        \]
        \begin{align*}
        \CN-Spec(K_{1, 13p^4-12p^3+11p^2-12p})&=\Big{\{}(0)^1, (-1)^{13p^4-12p^3+11p^2-12p-1},  \\
        & \qquad \qquad \qquad \qquad  (13p^4-12p^3+11p^2-12p-1)^1\Big{\}},
        \end{align*}
	Since $\CN-Spec(\mathcal{B}(G))$ is the union of the above sets, we get the required expression for $\CN-Spec(\mathcal{B}(G))$. 
    \end{proof}

    \begin{theorem}\label{all_energy_of_B(Q_4p)}
 Let $G=Q_{4p}$  and $p$ be any prime.
        \begin{enumerate}
            \item If $p=2$ then $E(\mathcal{B}(G))=2+6\sqrt{3}+4\sqrt{6}$,
	$LE(\mathcal{B}(G))= LE^+(\mathcal{B}(G))=\frac{4132}{35}$ and $E_{CN}(\mathcal{B}(G))=116$.
            \item If $p \geq 3$ then $E(\mathcal{B}(G))= 2+2\sqrt{3}+2p\sqrt{12}+2\sqrt{p^2-1}+2\sqrt{3p^2-3}+2\sqrt{12p^2-12p}$,
	$LE(\mathcal{B}(G))= LE^+(\mathcal{B}(G))=\frac{512p^4+2p^2+20p+50}{16p^2+p+5}$ 
	\quad  and
	$E_{CN}(\mathcal{B}(G))=32p^2-2p-10$.
        \end{enumerate}
    \end{theorem}
    \begin{proof}
(a) If $p=2$ then, by Theorem \ref{structure_of_Q_4p}, we have $\mathcal{B}(G)=K_2 \sqcup K_{1, 3} \sqcup 3K_{1, 12} \sqcup K_{1, 24}$. Now  by Lemma \ref{all_energy_of_star} we get
	\[
	E(K_2) = 2, \,  E( K_{1, 3}) = 2\sqrt{3}, \, E(3K_{1, 12})=2\sqrt{12} \, \text{ and } E(K_{1, 24})=2\sqrt{24}.
	\]
	Since $E(\mathcal{B}(G)) = E(K_2)+ E(K_{1, 3})+ E(3K_{1, 12})+ E(K_{1, 24})$ we get the required expression for $E(\mathcal{B}(G))$. By Theorem \ref{all_spec_of_B(Q_4p)}, we have $\L-Spec(\mathcal{B}(G))=\left\{(0)^6, (1)^{58}, (2)^1,\right.$ $\left. (4)^1, (13)^3, (25)^1 \right\}= \Q-Spec(\mathcal{B}(G))$, $m=64$ and $n=70$. Now, $\left|0-\frac{2(64)}{70}\right|=\frac{64}{35}$, $\left|1-\frac{2(64)}{70}\right|=\frac{29}{35}$, $\left|2-\frac{2(64)}{70}\right|=\frac{6}{35}$, $\left|4-\frac{2(64)}{70}\right|=\frac{76}{35}$. $\left|13-\frac{2(64)}{70}\right|=\frac{391)}{35}$ and $\left|25-\frac{2(64)}{70}\right|=\frac{811}{35}$. Therefore, by definition of $LE$ and $LE^+$ we have
 \begin{align*}
     LE(\mathcal{B}(G))=LE^+(\mathcal{B}(G))&=6 \times \frac{64}{35}+58 \times \frac{29}{35}+\frac{6}{35}+\frac{76}{35}+3 \times \frac{391}{35}+\frac{811}{35} \\
     &=\frac{4132}{35}.
 \end{align*}
 Finally,
\[
E_{CN}(K_2) = 0, \,  E_{CN}(K_{1, 3}) = 4, \, E_{CN}(3K_{1, 12}) =66 \text{ and }  E_{CN}(K_{1, 24}) = 46.	
\]
	Since $E_{CN}(\mathcal{B}(G)) = E_{CN}(K_2)+ E_{CN}(K_{1, 3})+ E_{CN}(3K_{1, 12})+ E_{CN}(K_{1, 24})$ we get the required expression for $E_{CN}(\mathcal{B}(G))$.
	
(b) If $p \geq 3$ then, by Lemma \ref{all_energy_of_star}, we get
	\[
	E(K_2) = 2, \,  E(K_{1, 3}) = 2\sqrt{3}, \, E(pK_{1, 12})=2p\sqrt{12}, \, E(K_{1, p^2-1}) = 2\sqrt{p^2 - 1}, 
        \]
        \[
         E(K_{1, 3p^2-3}) = 2\sqrt{3p^2 - 3} \text{ and }
         E(K_{1, 12p^2-12p})= 2\sqrt{12p^2-12p}.
         \]
Since $E(\mathcal{B}(G))\!=\!E(K_2)+ E(K_{1, 3})+ E(pK_{1, 12})+E(K_{1, p^2-1})+E(K_{1, 3p^2-3}) + E(K_{1, 12p^2-12p})$ 

\noindent we get the required expression for $E(\mathcal{B}(G))$. By Theorem \ref{all_spec_of_B(Q_4p)}, we have $\L-Spec(\mathcal{B}(G))=$ $\left\{\!(0)^{p+5}, (1)^{16p^2-p-5}, (2)^1, (4)^1, (13)^p, (p^2)^1, (3p^2-2)^1, (12p^2-12p+1)^1 \right\}\!=\!\Q-Spec(\mathcal{B}(G))$, 
$m=16p^2$ and $n=16p^2+p+5$. Now, $\left |0-\frac{2(16p^2)}{16p^2+p+5}\right|=\frac{32p^2}{16p^2+p+5}$, $\left |1-\frac{2(16p^2)}{16p^2+p+5}\right |=\frac{16p^2-p-5}{16p^2+p+5}$, $\left |2-\frac{2(16p^2)}{16p^2+p+5}\right|=\frac{2p+10}{16p^2+p+5}$, \quad $ \left |4-\frac{2(16p^2)}{16p^2+p+5}\right |=\frac{32p^2+4p+20}{16p^2+p+5}$, \quad $\left |13-\frac{2(16p^2)}{16p^2+p+5}\right |=\frac{176p^2+13p+65}{16p^2+p+5}$, 

\noindent $\left |p^2-\frac{2(16p^2)}{16p^2+p+5}\right| \quad = \, \frac{16p^4+p^3-27p^2}{16p^2+p+5}$, \quad $\left |(3p^2-2)-\frac{2(16p^2)}{16p^2+p+5}\right| \, = \quad \frac{48p^4+3p^3-49p^2-2p-10}{16p^2+p+5}$ \quad and 

\noindent $\left |(12p^2-12p+1)-\frac{2(16p^2)}{16p^2+p+5}\right | \quad = \quad\frac{192p^4-180p^3+32p^2-59p+5}{16p^2+p+5}$. Therefore, by definition of (signless) Laplacian energy we have 
 \begin{align*}
     LE(\mathcal{B}(G))&=LE^+(\mathcal{B}(G))\\
     &=(p+5)\frac{32p^2}{16p^2+p+5}+(16p^2-p-5)\frac{16p^2-p-5}{16p^2+p+5}+\frac{2p+10}{16p^2+p+5} \\
     & \qquad \qquad +\frac{32p^2+4p+20}{16p^2+p+5}+p\frac{176p^2+13p+65}{16p^2+p+5}+\frac{16p^4+p^3-27p^2}{16p^2+p+5} \\
     & \qquad +\frac{48p^4+3p^3-49p^2-2p-10}{16p^2+p+5}+\frac{192p^4-180p^3+32p^2-59p+5}{16p^2+p+5} \\
     &=\frac{512p^4+2p^2+20p+50}{16p^2+p+5}.
 \end{align*}
Finally,
\[
E_{CN}(K_2) = 0, \,  E_{CN}(K_{1, 3}) = 4, \, E_{CN}(pK_{1, 12}) = 22p, \,  E_{CN}(K_{1, p^2-1}) = 2p^2 - 4,
\]
\[
E_{CN}(K_{1, 3p^2-3}) = 6p^2-8 \, \text{ and }  E_{CN}(K_{1, 12p^2-12p}) = 24p^2-24p-2.	\]
	Since $E_{CN}(\mathcal{B}(G)) = E_{CN}(K_2)+ E_{CN}(K_{1, 3})+E_{CN}(pK_{1, 12})+ E_{CN}(K_{1, p^2-1})+ E_{CN}(K_{1, 3p^2-3})$ $+E_{CN}(K_{1, 12p^2-12p})$ we get the required expression for $E_{CN}(\mathcal{B}(G))$.
    \end{proof}
    \begin{theorem}\label{all_energy_of_B(Q_4p^2)}
Let $G=Q_{4p^2}$ and $p$ be any prime.
        \begin{enumerate}
            \item If $p=2$ then $E(\mathcal{B}(G)) =2+30\sqrt{3}+16\sqrt{6}$,
	$LE(\mathcal{B}(G))= LE^+(\mathcal{B}(G))=\frac{131314}{267}$ and $E_{CN}(\mathcal{B}(G))=490$.
            \item If $p \geq 3$ then $E(\mathcal{B}(G))= 2+2\sqrt{3}+2p^2\sqrt{12}+ 2\sqrt{p^2-1}+ 2\sqrt{3p^2-3}+2\sqrt{3p^4-3p^2}$ $+2(p-1)\sqrt{12p^2-12p}+2\sqrt{13p^4-12p^3+11p^2-12p}$,
	$LE(\mathcal{B}(G))= LE^+(\mathcal{B}(G))=\frac{512p^8+16p^5-24p^4-44p^3+118p^2-32p+54}{16p^4+p^2+p+5}$ 
	\quad  and
	$E_{CN}(\mathcal{B}(G))=32p^4-2p^2-2p-10$.
        \end{enumerate}
    \end{theorem}
    \begin{proof}
(a) If $p=2$ then, by Theorem \ref{structure_of_Q_4p^2}, we have $\mathcal{B}(G)=K_2 \sqcup K_{1, 3} \sqcup 5K_{1, 12} \sqcup 2K_{1, 24} \sqcup K_{1, 48} \sqcup K_{1, 96}$. Now  by Lemma \ref{all_energy_of_star} we get
	\[
	E(K_2) = 2, \,  E( K_{1, 3}) = 2\sqrt{3}, \, E(5K_{1, 12})=10\sqrt{12}, \, E(2K_{1, 24})=4\sqrt{24}, \, E( K_{1, 48}) = 2\sqrt{48}
        \]
        \[
        \text{ and } E( K_{1, 96}) = 2\sqrt{96}.
	\]
	Since $E(\mathcal{B}(G)) = E(K_2)+ E(K_{1, 3})+ E(5K_{1, 12})+ E(2K_{1, 24})+ E( K_{1, 48})+E( K_{1, 96})$ we get the required expression for $E(\mathcal{B}(G))$. By Theorem \ref{all_spec_of_B(Q_4p^2)}, we have $\L-Spec(\mathcal{B}(G))=\left\{(0)^{11}, (1)^{245}, (2)^1, (4)^1, (13)^5, (25)^2, (49)^1, (97)^1 \right\}= \Q-Spec(\mathcal{B}(G))$, $m=256$ and $n=267$. Now, $\left|0-\frac{2(256)}{267}\right|=\frac{512}{267}$, $\left|1-\frac{2(256)}{267}\right|=\frac{245}{267}$, $\left|2-\frac{2(256)}{267}\right|=\frac{22}{267}$, $\left|4-\frac{2(256)}{267}\right|=\frac{556}{267}$, $\left|13-\frac{2(256)}{267}\right|=\frac{2959)}{267}$, $\left|25-\frac{2(256)}{267}\right|=\frac{6163}{267}$, $\left|49-\frac{2(256)}{267}\right|=\frac{12571}{267}$ and $\left|97-\frac{2(256)}{267}\right|$ $=\frac{25387}{267}$. Therefore, by definition of $LE$ and $LE^+$ we have
 \begin{align*}
     LE(\mathcal{B}(G))=LE^+(\mathcal{B}(G))&=11 \times \frac{512}{267}+245 \times \frac{245}{267}+\frac{22}{267}+\frac{556}{267}+5 \times \frac{2959}{267} \\ 
     & \qquad \qquad \qquad \qquad +2 \times \frac{6163}{267}+ \frac{12571}{267}+\frac{25387}{267}\\
     &=\frac{131314}{267}.
 \end{align*}
 Finally,
\[
E_{CN}(K_2) = 0, \,  E_{CN}(K_{1, 3}) = 4, \, E_{CN}(5K_{1, 12}) =110, \, E_{CN}(2K_{1, 24}) =92, 
\]
\[
E_{CN}(K_{1, 48}) = 94 \, \text{ and } E_{CN}(K_{1, 96}) = 190.	
\]
	Since \, \, $E_{CN}(\mathcal{B}(G)) = E_{CN}(K_2)+ E_{CN}(K_{1, 3})+ E_{CN}(5K_{1, 12})+ E_{CN}(2K_{1, 24})$ $+E_{CN}(K_{1, 48})+E_{CN}(K_{1, 96})$ we get the required expression for $E_{CN}(\mathcal{B}(G))$.
	
(b) If $p \geq 3$ then, by Lemma \ref{all_energy_of_star}, we get
	\[
	E(K_2) = 2, \,  E(K_{1, 3}) = 2\sqrt{3}, \, E(p^2K_{1, 12})=2p^2\sqrt{12}, \, E(K_{1, p^2-1}) = 2\sqrt{p^2 - 1}, 
        \]
        \[
         E(K_{1, 3p^2-3}) = 2\sqrt{3p^2 - 3}, \,  E(K_{1, 3p^4-3^2}) = 2\sqrt{3p^4 - 3p^2}, 
         \]
         \[
         E((p-1)K_{1, 12p^2-12p})= 2(p-1)\sqrt{12p^2-12p} \,\text{ and }
         \]
         \[
          E(K_{1, 13p^4-12p^3+11p^2-12p}) = 2\sqrt{13p^4-12p^3+11p^2-12p}.
         \]
	Since $E(\mathcal{B}(G))$ is the sum of the above energies, we get the required expression for $E(\mathcal{B}(G))$. By Theorem \ref{all_spec_of_B(Q_4p^2)}, we have $\L-Spec(\mathcal{B}(G))=\Big{\{} (0)^{p^2+p+5}, (1)^{16p^4-p^2-p-5},  (2)^1, (4)^1$, $(13)^{p^2}, (p^2)^1, (3p^2-2)^1, (3p^4-3p^2+1)^1, (12p^2-12p+1)^{p-1},  (13p^4-12p^3+11p^2-12p+1)^1 \Big{\}}= \Q-Spec(\mathcal{B}(G))$, $m=16p^4$ and $n=16p^4+p^2+p+5$. Now, $\left |0-\frac{2(16p^4)}{16p^4+p^2+p+5}\right|=\frac{32p^4}{16p^4+p^2+p+5}$, 
	$\left |1-\frac{2(16p^4)}{16p^4+p^2+p+5}\right |=\frac{16p^4-p^2-p-5}{16p^4+p^2+p+5}$, 
	$\left |2-\frac{2(16p^4)}{16p^4+p^2+p+5}\right |=\frac{2p^2+2p+10}{16p^4+p^2+p+5},$ \quad 
	$\left |4-\frac{2(16p^4)}{16p^4+p^2+p+5}\right |=\frac{32p^4+4p^2+4p+20}{16p^4+p^2+p+5}$, 
	$\left |13-\frac{2(16p^4)}{16p^4+p^2+p+5}\right| = $ $\frac{176p^4+13p^2+13p+65}{16p^4+p^2+p+5}$, 
	$\left |p^2-\frac{2(16p^4)}{16p^4+p^2+p+5}\right|\!=\!
	\frac{16p^6-31p^4+p^3+5p^2}{16p^4+p^2+p+5}$,  
	
\noindent	$\left |(3p^2-2)-\frac{2(16p^4)}{16p^4+p^2+p+5}\right |=\frac{48p^6-61p^4+3p^3+13p^2-2p-10}{16p^4+p^2+p+5}$, \quad
$\left |(3p^4-3p^2+1)-\frac{2(16p^4)}{16p^4+p^2+p+5}\right|$ $=\frac{48p^8-45p^6+3p^5-4p^4-3p^3-14p^2+p+5}{16p^4+p^2+p+5}$, $\left |(12p^2-12p+1)-\frac{2(16p^4)}{16p^4+p^2+p+5}\right |= \frac{192p^6-192p^5-4p^4+p^2-11p+1}{16p^4+p^2+p+5}$  and  

$\left|(13p^4-12p^3+11p^2-12p+1)-\frac{2(16p^4)}{16p^4+p^2+p+5}\right|$

\qquad\qquad\qquad\qquad\qquad\qquad\qquad\qquad\qquad $=\frac{208p^8-192p^7+189p^6-175p^5+48p^4-61p^3+44p^2-59p+5}{16p^4+p^2+p+5}$. 

\noindent Therefore, by definition of (signless) Laplacian energy we have 
 \begin{align*}
     & \qquad LE(\mathcal{B}(G))=LE^+(\mathcal{B}(G)) \\
     &=(p^2+p+5)\frac{32p^4}{16p^4+p^2+p+5}+\frac{(16p^4-p^2-p-5)^2}{16p^4+p^2+p+5}+\frac{2p^2+2p+10}{16p^4+p^2+p+5} \\
     & \qquad +\frac{32p^4+4p^2+4p+20}{16p^4+p^2+p+5}+p^2\frac{176p^4+13p^2+13p+65}{16p^4+p^2+p+5}+\frac{16p^6-31p^4+p^3+5p^2}{16p^4+p^2+p+5} \\
     & + \frac{48p^6-61p^4+3p^3+13p^2-2p-10}{16p^4+p^2+p+5} +\frac{192p^6-192p^5-4p^4+p^2-11p+1}{16p^4+p^2+p+5} \\
     & +\frac{48p^8-45p^6+3p^5-4p^4-3p^3-14p^2+p+5}{16p^4+p^2+p+5} \\
     & +\frac{208p^8-192p^7+189p^6-175p^5+48p^4-61p^3+44p^2-59p+5}{16p^4+p^2+p+5} \\
     &=\frac{512p^8+16p^5-24p^4-44p^3+118p^2-32p+54}{16p^4+p^2+p+5}.
 \end{align*}
Finally,
\[
E_{CN}(K_2) = 0, \,  E_{CN}(K_{1, 3}) = 4, \, E_{CN}(p^2K_{1, 12}) = 22p^2, \,  E_{CN}(K_{1, p^2-1}) = 2p^2 - 4,
\]
\[
E_{CN}(K_{1, 3p^2-3}) = 6p^2-8, \, E_{CN}(K_{1, 3p^4-3p^2}) = 6p^4-6p^2-2, 
\]
\[
E_{CN}((p-1)K_{1, 12p^2-12p}) = 24p^3-48p^2+22p+2 \text{ and }
\]
\[
E_{CN}(K_{1, 13p^4-12p^3+11p^2-12p}) = 26p^4-24p^3+22p^2-24p-2.
\]
	Since $E_{CN}(\mathcal{B}(G))\!=\!E_{CN}(K_2)+ E_{CN}(K_{1, 3})+E_{CN}(p^2K_{1, 12})+ E_{CN}(K_{1, p^2-1})+ E_{CN}(K_{1, 3p^2-3})$ $+E_{CN}(K_{1, 3p^4-3p^2})+E_{CN}(K_{1, 12p^2-12p})+E_{CN}(K_{1, 13p^4-12p^3+11p^2-12p})$ we get the required expression for $E_{CN}(\mathcal{B}(G))$.
This completes the proof.
    \end{proof}
    \begin{theorem}\label{all_energy_comparison_Q_4p}
	If $G=Q_{4p}$ then $\mathcal{B}(G)$ is hypoenergetic but not hyperenergetic,  L-hyper-energetic, Q-hyperenergetic and CN-hyperenergetic.
\end{theorem}
\begin{proof}
	\textbf{Case 1.} $p=2$
	
	 By Theorem \ref{structure_of_Q_4p} and Theorem \ref{all_energy_of_B(Q_4p)}, we have $|V(\mathcal{B}(G))|=70$ and $E(\mathcal{B}(G)) \approx 36.0466$. Clearly, $|V(\mathcal{B}(G))|> E(\mathcal{B}(G))$.	Thus, $\mathcal{B}(G)$ is hypoenergetic. 
	
	We have $E(K_{|V(\mathcal{B}(G))|}) = E(K_{70})= 2(70-1)=138 >36.0466 \approx E(\mathcal{B}(G))$. Therefore,  $\mathcal{B}(G)$ is not hyperenergetic.

	Also, $LE(K_{70})=LE^+(K_{70})=138$. From Theorem \ref{all_energy_of_B(Q_4p)}, we have $LE(\mathcal{B}(G))=LE^+(\mathcal{B}(G)) \approx 118.0571$.
Therefore, $LE^+(K_{70}) = LE(K_{70}) > LE(\mathcal{B}(G))=LE^+(\mathcal{B}(G))$.	Hence,  $\mathcal{B}(G)$ is neither L-hyperenergetic nor  Q-hyperenergetic.
	
 We have  $E_{CN}(K_{70})=2(70-1)(70-2)=9384 > 116 = E_{CN}(\mathcal{B}(G))$ (using Theorem \ref{all_energy_of_B(Q_4p)}). Hence, $\mathcal{B}(G)$ is not CN-hyperenergetic. 
 
  \noindent \textbf{Case 2.}  $p \geq 3$
  
  By Theorem \ref{structure_of_Q_4p} and Theorem \ref{all_energy_of_B(Q_4p)}, we have $|V(\mathcal{B}(G))|=16p^2+p+5$ and
	\[
	E(\mathcal{B}(G))=2+2\sqrt{3}+2p\sqrt{12}+2\sqrt{p^2-1}+2\sqrt{3p^2-3}+2\sqrt{12p^2-12p}.
	\]
	Since $4>2\sqrt{3}, \, 13p> 2p\sqrt{12}, \, p^2>2\sqrt{p^2-1}, \, 3p^2-2>2\sqrt{3p^2-3}$ and $12p^2-12p+1 >2\sqrt{12p^2-12p}$ \,\, we have 
	\begin{align*}
	    4+13p+p^2+3p^2-2+12p^2-12p+1+2 > 2+2\sqrt{3}+2p\sqrt{12}&+2\sqrt{p^2-1}+2\sqrt{3p^2-3} \\ 
        &+2\sqrt{12p^2-12p}.
	\end{align*}
	Therefore,
	\begin{equation} \label{Q-4p-hypo}
		|V(\mathcal{B}(G))|> E(\mathcal{B}(G)).
	\end{equation}	 
	Thus, $\mathcal{B}(G)$ is hypoenergetic. 
	
	We have $E(K_{|V(\mathcal{B}(G))|}) = E(K_{16p^2+p+5})= 2(16p^2+p+5-1)=32p^2+2p+8 >16p^2+p+5 = |V(\mathcal{B}(G))|> E(\mathcal{B}(G))$ (using \eqref{Q-4p-hypo}). Therefore,  $\mathcal{B}(G)$ is not hyperenergetic.

	Also, $LE(K_{16p^2+p+5})=LE^+(K_{16p^2+p+5})=32p^2+2p+8$. From Theorem \ref{all_energy_of_B(Q_4p)}, we have $LE(\mathcal{B}(G))=LE^+(\mathcal{B}(G))=\frac{512p^4+2p^2+20p+50}{16p^2+p+5}$. Now, 
	\begin{align*}
		LE(K_{16p^2+p+5})-LE(\mathcal{B}(G))&= LE^+(K_{16p^2+p+5})-LE^+(\mathcal{B}(G))\\ &=32p^2+2p+8-\frac{512p^4+2p^2+20p+50}{16p^2+p+5} \\
		&= \frac{64p^3+288p^2-2p-10}{16p^2+p+5} > 0.
	\end{align*}
Therefore, 
	\[
	LE^+(K_{16p^2+p+5}) = LE(K_{16p^2+p+5}) > LE(\mathcal{B}(G))=LE^+(\mathcal{B}(G)).
	\] 
	Hence,  $\mathcal{B}(G)$ is neither L-hyperenergetic nor  Q-hyperenergetic.
	
	We have  $E_{CN}(K_{16p^2+p+5})=2(16p^2+p+5-1)(16p^2+p+5-2)=(32p^2+2p+8)(16p^2+p+3) > 32p^2-2p-10 = E_{CN}(\mathcal{B}(G))$ (using Theorem \ref{all_energy_of_B(Q_4p)}). Hence, $\mathcal{B}(G)$ is not CN-hyperenergetic. 
 
  This completes the proof.
\end{proof}
\begin{theorem}\label{all_energy_comparison_Q_4p^2}
	If $G=Q_{4p^2}$, then $\mathcal{B}(G)$ is hypoenergetic but not hyperenergetic,  L-hyperenergetic, Q-hyperenergetic and CN-hyperenergetic.
\end{theorem}
\begin{proof}
	\textbf{Case 1.} $p=2$
	
	By Theorem \ref{structure_of_Q_4p^2} and Theorem \ref{all_energy_of_B(Q_4p)}, we have $|V(\mathcal{B}(G))|=267$ and $E(\mathcal{B}(G)) \approx 93.1533$. Clearly, $|V(\mathcal{B}(G))|> E(\mathcal{B}(G))$.	Thus, $\mathcal{B}(G)$ is hypoenergetic. 
	
	We have $E(K_{|V(\mathcal{B}(G))|}) = E(K_{267})= 2(267-1)=532 >93.1533 \approx E(\mathcal{B}(G))$. Therefore,  $\mathcal{B}(G)$ is not hyperenergetic.
		
	Also, $LE(K_{267})=LE^+(K_{267})=532$. From Theorem \ref{all_energy_of_B(Q_4p^2)}, we have $LE(\mathcal{B}(G))$ $=LE^+(\mathcal{B}(G)) \approx 491.8127$. Therefore, $LE^+(K_{267})\!=\!LE(K_{267}) > LE(\mathcal{B}(G))=LE^+(\mathcal{B}(G))$.	Hence,  $\mathcal{B}(G)$ is neither L-hyperenergetic nor  Q-hyperenergetic.
	
 We have  $E_{CN}(K_{267})=2(267-1)(267-2)=140980 > 490 = E_{CN}(\mathcal{B}(G))$ (using Theorem \ref{all_energy_of_B(Q_4p^2)}). Hence, $\mathcal{B}(G)$ is not CN-hyperenergetic. 
 
  \noindent \textbf{Case 2.} $p \geq 3$
  
  By Theorem \ref{structure_of_Q_4p^2} and Theorem \ref{all_energy_of_B(Q_4p^2)}, we have $|V(\mathcal{B}(G))|=16p^4+p^2+p+5$ and
	\begin{align*}
	E(\mathcal{B}(G))=2+2\sqrt{3}&+2p^2\sqrt{12}+2\sqrt{p^2-1}+2\sqrt{3p^2-3}+2\sqrt{3p^4-3p^2} \\
 & +2(p-1)\sqrt{12p^2-12p}+2\sqrt{13p^4-12p^3+11p^2-12p}.
	\end{align*}
	Since $4>2\sqrt{3}, \, 13p^2> 2p^2\sqrt{12}, \, p^2>2\sqrt{p^2-1}, \, 3p^2-2>2\sqrt{3p^2-3}, \, 3p^4-3p^2+1>2\sqrt{3p^4-3p^2}, \, (p-1)(12p^2-12p+1) >2(p-1)\sqrt{12p^2-12p}$ and $13p^4-12p^3+11p^2-12p+1>2\sqrt{13p^4-12p^3+11p^2-12p}$, by adding both sides we have
	\begin{equation} \label{Q-4p^2-hypo}
		|V(\mathcal{B}(G))|> E(\mathcal{B}(G)).
	\end{equation}	 
	Thus, $\mathcal{B}(G)$ is hypoenergetic. 
	
	We have $E(K_{|V(\mathcal{B}(G))|}) = E(K_{16p^4+p^2+p+5})= 2(16p^4+p^2+p+5-1)=32p^4+2p^2+2p+8 >16p^4+p^2+p+5 = |V(\mathcal{B}(G))|> E(\mathcal{B}(G))$ (using \eqref{Q-4p^2-hypo}). Therefore,  $\mathcal{B}(G)$ is not hyperenergetic.

	Also, $LE(K_{16p^4+p^2+p+5})=LE^+(K_{16p^4+p^2+p+5})=32p^4+2p^2+2p+8$. From Theorem \ref{all_energy_of_B(Q_4p^2)}, we have $LE(\mathcal{B}(G))=LE^+(\mathcal{B}(G))=\frac{512p^8+16p^5-24p^4-44p^3+118p^2-32p+54}{16p^4+p^2+p+5}$. Now, 
	\begin{align*}
		LE(&K_{16p^4+p^2+p+5})-LE(\mathcal{B}(G))= LE^+(K_{16p^4+p^2+p+5})-LE^+(\mathcal{B}(G))\\ &=32p^4+2p^2+2p+8-\frac{512p^8+16p^5-24p^4-44p^3+118p^2-32p+54}{16p^4+p^2+p+5} \\
		&= \frac{64p^6+48p^5+314p^4+48p^3-98p^2+50p-14 }{16p^4+p^2+p+5} > 0.
	\end{align*}
Therefore, 
	\[
	LE^+(K_{16p^4+p^2+p+5}) = LE(K_{16p^4+p^2+p+5}) > LE(\mathcal{B}(G))=LE^+(\mathcal{B}(G)).
	\] 
	Hence,  $\mathcal{B}(G)$ is neither L-hyperenergetic nor  Q-hyperenergetic.
	
	We have  $E_{CN}(K_{16p^4+p^2+p+5})=2(16p^4+p^2+p+5-1)(16p^4+p^2+p+5-2)=(32p^4+2p^2+2p+8)(16p^4+p^2+p+3) > 32p^4-2p^2-2p-10 = E_{CN}(\mathcal{B}(G))$ (using Theorem \ref{all_energy_of_B(Q_4p^2)}). Hence, $\mathcal{B}(G)$ is not CN-hyperenergetic. 
 
  This completes the proof.
\end{proof}
\begin{theorem}
	If $G=Q_{4p}$, then $E(\mathcal{B}(G))<LE(\mathcal{B}(G))=LE^+(\mathcal{B}(G))$.
\end{theorem}
\begin{proof}
    \textbf{Case 1.} $p=2$.

     By Theorem \ref{all_energy_of_B(Q_4p)} (a), we get $E(\mathcal{B}(G)) \approx 36.0466$ and $LE(\mathcal{B}(G))= LE^+(\mathcal{B}(G)) \approx 118.0571$. Clearly, $LE(\mathcal{B}(G))= LE^+(\mathcal{B}(G))>E(\mathcal{B}(G))$.
     
    \noindent \textbf{Case 2.} $p \geq 3$.

     By Theorem \ref{all_energy_of_B(Q_4p)} (b), we get
	\[
	E(\mathcal{B}(G))= 2+2\sqrt{3}+2p\sqrt{12}+2\sqrt{p^2-1}+2\sqrt{3p^2-3}+2\sqrt{12p^2-12p} \qquad \text{ and }
        \]
	\[
        LE(\mathcal{B}(G))= LE^+(\mathcal{B}(G))=\frac{512p^4+2p^2+20p+50}{16p^2+p+5}.
	\] 
	
 Also, from Theorem \ref{all_energy_comparison_Q_4p}, we have $|V(\mathcal{B}(G))|=16p^2+p+5> E(\mathcal{B}(G))$. Now, 
 \begin{align*}
     LE(\mathcal{B}(G))-|V(\mathcal{B}(G))|&= LE^+(\mathcal{B}(G))-|V(\mathcal{B}(G))| \\
     &=\frac{512p^4+2p^2+20p+50}{16p^2+p+5}-16p^2+p+5 \\
     &= \frac{256p^4-32p^3-159p^2+10p+25}{16p^2+p+5} \\
     &= \frac{p^2(256p^2-32p-159)+10p+25}{16p^2+p+5}:=\frac{f(p)}{g(p)}.
 \end{align*}
 Since, for all prime $p$, $256p^2-32p>159$  and $16p^2+p+5>0$ we have $\frac{f(p)}{g(p)}>0$. Hence, $LE(\mathcal{B}(G))= LE^+(\mathcal{B}(G))>|V(\mathcal{B}(G))|>E(\mathcal{B}(G))$.
\end{proof}
\begin{theorem}
	If $G=Q_{4p^2}$, then $E(\mathcal{B}(G))<LE(\mathcal{B}(G))=LE^+(\mathcal{B}(G))$.
\end{theorem}
\begin{proof}
    \textbf{Case 1.} $p=2$.

     By Theorem \ref{all_energy_of_B(Q_4p^2)} (a), we get $E(\mathcal{B}(G)) \approx 93.1533$ and $LE(\mathcal{B}(G))= LE^+(\mathcal{B}(G)) \approx 491.8127$. Clearly, $LE(\mathcal{B}(G))= LE^+(\mathcal{B}(G))>E(\mathcal{B}(G))$.
     
    \noindent \textbf{Case 2.} $p \geq 3$.

     By Theorem \ref{all_energy_of_B(Q_4p^2)} (b), we get
	\begin{align*}
	E(\mathcal{B}(G))= 2+2\sqrt{3}&+2p^2\sqrt{12}+2\sqrt{p^2-1}+2\sqrt{3p^2-3}+2\sqrt{3p^4-3p^2} \\
 & +2(p-1)\sqrt{12p^2-12p}+2\sqrt{13p^4-12p^3+11p^2-12p} \qquad \text{ and }
        \end{align*}
	\[
        LE(\mathcal{B}(G))= LE^+(\mathcal{B}(G))=\frac{512p^8+16p^5-24p^4-44p^3+118p^2-32p+54}{16p^4+p^2+p+5}.
	\] 
	
 Also, from Theorem \ref{all_energy_comparison_Q_4p^2}, we have $|V(\mathcal{B}(G))|=16p^4+p^2+p+5> E(\mathcal{B}(G))$. Now, 
 \begin{align*}
     LE(\mathcal{B}(G))&-|V(\mathcal{B}(G))|= LE^+(\mathcal{B}(G))-|V(\mathcal{B}(G))| \\
     &=\frac{512p^8+16p^5-24p^4-44p^3+118p^2-32p+54}{16p^4+p^2+p+5}-16p^4+p^2+p+5 \\
     &= \frac{256p^8-32p^6-16p^5-185p^4-46p^3+107p^2-42p+29}{16p^4+p^2+p+5} \\
     &= \frac{p^3(256p^5-32p^3-16p^2-185p-46)+107p^2-42p+29}{16p^4+p^2+p+5}:=\frac{f(p)}{g(p)}.
 \end{align*}
 Since, for all prime $p$, $256p^5-32p^3-16p^2-185p>46$  and $16p^4+p^2+p+5>0$ we have $\frac{f(p)}{g(p)}>0$. Hence, $LE(\mathcal{B}(G))= LE^+(\mathcal{B}(G))>|V(\mathcal{B}(G))|>E(\mathcal{B}(G))$.
\end{proof}
\section*{Concluding remark}
We conclude this paper with the following remark.
\begin{remark}\label{last-rem}
	If $G = D_{2p}, D_{2p^2}, Q_{4p}$ and $Q_{4p^2}$, where $p$ is any prime, then
	\begin{enumerate}
		\item 	 $\mathcal{B}(G)$  is not integral but L-integral, Q-integral and CN-integral.
		\item $\mathcal{B}(G)$ is hypoenergetic but not hyperenergetic,  L-hyperenergetic, Q-hyperenergetic and CN-hyperenergetic.
		\item $\mathcal{B}(G)$ satisfies E-LE conjecture. 
	\end{enumerate}
\end{remark}
It may be interesting to consider the following problem.
\begin{prob}
	Is there any finite non-abelian group $G$ such that $\mathcal{B}(G)$	does not hold the conclusions of Remark \ref{last-rem}?
\end{prob}

\vspace{.5cm}

\noindent {\bf Acknowledgement.} The first author is thankful to Council of Scientific and Industrial Research  for the fellowship (File No. 09/0796(16521)/2023-EMR-I).

\section*{Declarations}
\begin{itemize}
	\item Funding: No funding was received by the authors.
	\item Conflict of interest: The authors declare that they have no conflict of interest.
	\item Availability of data and materials: No data was used in the preparation of this manuscript.
\end{itemize}


\begin{thebibliography}{10}

\bibitem{AEN17}
Abdussakir, --, Elvierayani, R. R. and Nafisah, M.  On the spectra of commuting and non-commuting graph on dihedral group, \emph{Cauchy-Jurnal Matematika Murni dan Aplikasi}, \textbf{4}(4), 176--182, 2017.

 \bibitem{ACGMR11} 
Abreu, N., Cardoso, D. M., Gutman, I., Martins, E. A. and Robbiano, M. Bounds for the signless Laplacian energy. \emph{Linear Algebra and its Applications}, \textbf{435}(10), 2365--2374, 2011.

\bibitem{AL2021}
Ali, F.,  and Li, Y. The connectivity and the spectral radius of commuting on certain groups, {\em Linear and Multilinear Algebra}, \textbf{69}(16), 2945--2958, 2021.

\bibitem{AAGS12}
Alwardi, A., Arsic, B., Gutman, I., Sonera, N.D. The common neighborhood graph and its energy, {\em Iranian Journal of Mathematical Sciences and Informatics}, \textbf{7}(2), 1--8, 2012.

\bibitem{ASG}
Alwardi, A.,   Soner, N. D. and   Gutman, I.  On the common-neighborhood energy of a graph, 
{\em Bulletin (Acad$\acute{\rm e}$mie Serbe Des Sciences Et Des Arts. Classe Des Sciences Math$\acute{\rm e}$matiques Et Naturelles. Sciences Math$\acute{\rm e}$matiques)}, 
{\bf 36}, 49--59, 2011.

\bibitem{ANC-2022}
Arunkumar, G.,  Cameron, P. J.,  Nath, R. K. and  Selvaganesh, L. Super graphs on groups I, {\em Graphs and Combinatorics}, \textbf{38}, Article No. 100, 2022.

\bibitem{ACGN24}
Arunkumar, G.,   Cameron, P. J.,   Ganeshbabu, R. and  Nath, R. K. Main functions and the spectrum of super graphs, preprint, https://doi.org/10.48550/arXiv.2408.00390

\bibitem{BaCv-2002}
Balińska, K.,  Cvetković, D.,  Radosavljević, Z.,  Simić, S. and  Stevanović, D. A survey on integral graphs, {\em Univ. Beograd. Publ. Elektrotehn. Fak. Ser. Mat.} \textbf{13}, 42--65, 2002.

\bibitem{BA23}
Banerjee, S. and   Adhikari, A. On spectra of power graphs of finite cyclic and dihedral groups, {\em Rocky Mountain Journal of  Mathematics}, \textbf{53}(2), 341--356, 2023. 


\bibitem{BN21}
Bhowal, P. and  Nath, R. K. Spectral aspects of commuting conjugacy class graph of finite groups, {\em Algebraic Structures and Their Applications}, \textbf{8}(2),  67--118, 2021.


\bibitem{BN24}
Bhowal, P. and  Nath, R. K. Spectrum and energies of commuting conjugacy class graph of a finite group, {\em  Algebraic Structures and Their Applications}, \textbf{12}(1),  33--49, 2025.

\bibitem{cameron2021graphs}
Cameron, P. J. Graphs defined on groups, {\em International Journal of Group Theory}, \textbf{11}(2), 53--107, 2022.

\bibitem{CaJa-2024}
Cameron, P. J.,   Jannat, F. E.,   Nath, R. K. and Sharafdini,  R. A survey on conjugacy class graphs of groups, {\em Expositiones Mathematicae}, \textbf{42}(4), Article No. 125585, 2024.	





\bibitem{CPA18}
Chattopadhyay, S.,  Panigrahi, P. and  Atik, F. Spectral radius of power graphs on certain finite groups, {\em Indagationes Mathematicae}, \textbf{29}(2),  730--737, 2018.

\bibitem{CRS-2007}
Cvetković, D.,  Rowlinson, P. and  Simić, S. K. Signless Laplacians of finite graphs, {\em Linear Algebra and its Application}, \textbf{423},  155--171, 2007.

\bibitem{DMP24}
Dalal, S., Mukherjee, S. and Patra, K. L. Spectrum of super commuting graphs of some finite groups, {\em Computational and Applied Mathematics}, \textbf{43}, Article No. 348, 2024. 




\bibitem{DEN-23}
Das,  S., Erfanian, A. and  Nath, R. K. On a bipartite graph defined on groups, {\em Journal of Algebra and Its Application}, https://doi.org/10.1142/S0219498826501926.



\bibitem{DEN-24}
Das,  S., Erfanian, A. and  Nath, R. K. Zagreb indices of the subgroup generating bipartite graph, {\em Discrete Mathematics, Algorithms and Applications}, https://doi.org/10.1142/S1793830925501411.

\bibitem{DBN-2020}
Dutta, P.,   Bagchi, B. and  Nath, R. K. Various energies of commuting graphs of finite nonabelian groups, {\em Khayyam Journal of Mathematics}, \textbf{6}(1), 27--45, 2020.

\bibitem{DDN18} 
Dutta, P., Dutta, J. and Nath, R. K. Laplacian spectrum of noncommuting graphs of finite groups,  \emph{Indian Journal of Pure and Applied Mathematics}, \textbf{49}(2), 205--216, 2018.

\bibitem{DN-2017}
Dutta, J. and  Nath, R. K. Spectrum of commuting graphs of some classes of finite groups, {\em Matematika}, \textbf{33}(1), 87--95, 2017.

\bibitem{DN-2017MV}
Dutta, J. and  Nath, R. K. Finite groups whose commuting graphs are integral, {\em Matematički Vesnik}, \textbf{69}(3), 226--230, 2017.

\bibitem{DN-2018}
Dutta, J. and  Nath, R. K. Laplacian and signless Laplacian spectrum of commuting graphs of finite groups, {\em Khayyam Journal of Mathematics}, \textbf{4}(1), 77--87, 2018.

\bibitem{DN-2021}
Dutta, J. and  Nath, R. K. Various energies of commuting graphs of some super integral groups, {\em Indian Journal of Pure and Applied Mathematics}, \textbf{52}(1), 1--10, 2021.

\bibitem{DN18}
Dutta, P. and Nath, R. K. On Laplacian energy of non-commuting graphs of finite groups,  \emph{Journal of Linear and Topological Algebra}, \textbf{7}(2), 121--132, 2018.



\bibitem{FKMN-2005}
Fallat, S.,   Kirkland, S.,   Molitierno, J. and  Neumann, M. On graphs whose Laplacian matrices have distinct integer eigenvalues, {\em Journal of Graph Theory} \textbf{50},  162--174, 2005.

\bibitem{FN24}
Fasfous, W. N. T. and  Nath, R. K. Inequalities involving energy and Laplacian energy of non-commuting graphs of finite groups, {\em Indian Journal of Pure and Applied Mathematics}, \textbf{56}(2), 791--812, 2025.


\bibitem{FSN-2021}
Fasfous, W. N. T.,   Sharafdini, R. and  Nath, R. K. Common neighborhood spectrum of commuting graphs of finite groups, {\em Algebra and Discrete Mathematics}, \textbf{32}(1),  33--48, 2021.

\bibitem{GGA17}
Ghorbani, M. and Gharavi-Alkhansari, Z. On the energy of non-commuting graphs,  \emph{Journal of Linear and Topological Algebra}, \textbf{6}(2), 135--146, 2017.

\bibitem{GGAZ17}
Ghorbani, M., Gharavi-Alkhansari, Z. and Bashi, A. Z. On the eigenvalues of non-commuting graphs,  \emph{Journal of Algebraic Structures and Their Applications}, \textbf{4}(2), 27--38, 2017.

\bibitem{Gutman-78}
Gutman, I. The energy of a graph, {\em Ber. Math-Statist. Sekt. Forschungsz. Graz}, \textbf{103}, 1--22, 1978.

\bibitem{Gutman-99}
Gutman, I.
\newblock Hyperenergetic molecular graphs,
\newblock {\em Journal of the Serbian Chemical Society}, \textbf{64}, 199--205, 1999.

\bibitem{E-LE-Gutman}
Gutman, I., Abreu, N. M. M., Vinagre, C. T. M., Bonifacioa, A. S. and Radenkovic, S. Relation between energy and Laplacian energy, {\em MATCH Communications in Mathematical and in Computer Chemistry}, \textbf{59}, 343--354, 2008.

\bibitem{GZ06} 
Gutman, I. and Zhou, B. Laplacian energy of a graph. \emph{Linear Algebra and its Applications}, \textbf{414}(1), 29--37, 2006.

\bibitem{HA17}
Hamzeh A. and Ashrafi, A. R.  Spectrum and L-spectrum of the power graph and its main supergraph for certain finite groups, {\em Filomat},  \textbf{31}(16),  5323--5334, 2017.

\bibitem{fns20}
Fasfous, W. N. T., Nath, R. K. and Sharafdini, R. Various spectra and  energies of commuting graphs of finite rings, 
\emph{Hacettepe Journal of Mathematics and Statistics}, \textbf{49}(6), 1915--1925, 2020.

\bibitem{hansen2007comparing}
Hansen, P. and Vuki{\v{c}}evi{\'c}, D. Comparing the Zagreb indices, {\em Croatica Chemica Acta}, \textbf{80}(2), 165--168, 2007.


\bibitem{HaSc-1974}
Harary, F. and  Schwenk, A. J. Which graphs have integral spectra?, {\em Graphs and Combinatorics} (R. Bari and F. Harary, eds.), Springer- Verlag, Berlin, \textbf{45}, 1974.





\bibitem{JN25}
Jannat, F. E. and  Nath, R. K. Common neighborhood spectrum and energy of commuting conjugacy class graph, {\em Journal of Algebraic Systems}, \textbf{12}(2),  301--326, 2025.

\bibitem{KSCC-2021}
Kumar, A., Selvaganesh, L., Cameron, P. J., and Chelvam, T. T. (2021). Recent developments on the power graph of finite groups - a survey, {\em AKCE International Journal of Graphs and Combinatorics}, \textbf{18}(2), 65--94, 2021. 

\bibitem{MEM-2023}
Mahtabi, M., Erfanian, A. and  Mahtabi, R. A bipartite graph associated to a finite group, {\em Asian-European Journal of Mathematics}, \textbf{16}(10), Article No. 2350185, 2023.

\bibitem{MGA17}
Mehranian, Z., Gholami, A., and Ashrafi, A. R.  The Spectra of power graphs of certain finite groups, {\em Linear and Multilinear Algebra}, \textbf{65}(5), 1003--1010, 2017. 

\bibitem{Merris-1994}
Merris, R. Degree maximal graphs are Laplacian integral, {\em Linear Algebra and its Application}, \textbf{199}, 381--389, 1994.


\bibitem{N-2018}
Nath, R. K. Various spectra of commuting graphs of $n$-centralizer finite groups, {\em International Journal of Engineering, Science and Technology}, \textbf{10}(2S), 170--172, 2018.

\bibitem{NFDS-2021}
Nath, R. K.,  Fasfous, W. N. T.,   Das, K. C. and  Shang, Y. Common neighborhood energy of commuting graphs of finite groups, {\em Symmetry}, \textbf{13}(9), Article No. 1651, 2021.

\bibitem{P19}
Panda, R. P. Laplacian spectra of power graphs of certain finite groups, {\em Graphs and Combinatorics} \textbf{35}, 1209--1223, 2019. 

\bibitem{RPCA23}
Rather, B. A.,   Pirzada, S.,  Chishti, T. A. and Alghamdi, A. M. On normalized Laplacian eigenvalues of power graphs associated to finite cyclic groups, {\em  Discrete Mathematics, Algorithms and Applications},  \textbf{15}(2), Article No. 2250070, 2023.

\bibitem{S22} 
Singh, S. N. Laplacian spectra of power graphs of certain prime-power Abelian groups, {\em  Asian-European Journal of Mathematics}, \textbf{15}(2), Article No. 2250026, 2022.


\bibitem{DN-2022}
Sharafdini, R.,   Nath, R. K. and  Darbandi, R. Energy of commuting graph of finite AC-groups, {\em Proyecciones Journal of Mathematics}, \textbf{41}(1), 263-273, 2022.

\bibitem{SN24} 
Sharma, M. and  Nath, R. K. Signless Laplacian energies of non-commuting graphs of finite groups and related results, {\em Discrete Mathematics, Algorithms and Applications}, \textbf{17}(4), Article No. 2450060, 2025.

\bibitem{Stanic-2009}
Stanić, Z. Some results on Q-integral graphs, {\em Ars Combin.} \textbf{90}, 321--335, 2009. 

\bibitem{TA19}
Torktaz, M. and  Ashrafi, A.R. Spectral properties of the commuting graphs of certain groups, {\em AKCE International Journal of Graphs and Combinatorics}, \textbf{16}(3), 300--309, 2019.

 \bibitem{Walikar-99}
Walikar, H. B., Ramane, H. S. and Hampiholi, P. R.
\newblock On the energy of a graph,
\newblock  {\em Graph Connections, Allied Publishers, New Delhi}, 120--123, 1999. 


\end{thebibliography}
\end{document}